\documentclass[reqno]{amsart}

\usepackage{amsmath}
\numberwithin{equation}{section}

\usepackage{amssymb}
\usepackage{amsthm}

\theoremstyle{plain}  % default setting: bold heading and italic body
\newtheorem{thm}{Theorem}[section]

\newtheorem{cor}[thm]{Corollary}
\newtheorem{lem}[thm]{Lemma}
\newtheorem{prop}[thm]{Proposition}

\theoremstyle{definition} % bold heading and normal body
\newtheorem{defn}{Definition}[section]

\theoremstyle{remark} % italic heading and normal body

\begin{document}

\title{An expression for the Homflypt polynomial and some applications}

\author{David Emmes}
\address{Unaffiliated}
\email{davidemmes@gmail.com}


\begin{abstract}
%% Text of abstract
Associated with each oriented link is the two variable Homflypt polynomial.
The Morton-Franks-Williams (MFW) inequality gives rise to an
expression for the Homflypt polynomial with MFW coefficient polynomials.
These MFW coefficient polynomials
are labelled in a braid-dependent manner and may be zero, but display
a number of interesting relations.
One consequence is an
expression for the first three Laurent coefficient polynomials in z as a function of 
the other coefficient polynomials and three
link invariants: the minimum v-degree and v-span of the Homflypt polynomial, and
the Conway polynomial. 

These expressions are used to derive additional properties of the Homflypt 
polynomial for general n-braid links. One specific result is that
the Jones and Homflypt polynomials distinguish the same three-braid links.
\end{abstract}

\maketitle

\section{Introduction}
The focus of this paper is to develop and display a relationship among the coefficient polynomials
of the \mbox{Homflypt} polynomial, and to use this to 
establish some properties of the \mbox{Homflypt} polynomial
for general $n$-braid links and three-braid links.
The background required for this paper is a very basic understanding of skein relations, 
skein polynomials, and braid groups. 
Such material may be found in any introductory textbook, as for example K.~Murasugi's,
\cite{71}, and may 
also be found in many excellent surveys, \cite{10}.

The expression for the Homflypt polynomial is based on the
Morton-Franks-Williams (MFW) inequality, \cite{30}, \cite{60}, that
provides lower and upper bounds for the allowed powers of $v$ in the
\mbox{Homflypt} polynomial. This gives rise to an indexed set of 
MFW polynomials whose properties form the basis for the paper.
Section~\ref{SubPropsSkeinPoly} defines these braid-dependent MFW polynomials.

The second section contains the primary results.
Thm.~\ref{Hcoeffrel} expresses the first three
MFW polynomials in terms of braid properties and the remaining MFW polynomials. 
This leads to Thm.~\ref{HLaurentcoeffrel} that establishes a system of relations
among MFW polynomials.
Cor.~\ref{HProperLcoeffrel} expresses the 
first three Laurent coefficient polynomials in terms of
the other coefficient polynomials and three
link invariants: the minimum $v$-degree and $v$-span of the \mbox{Homflypt} polynomial, and
the Conway polynomial.
Thm.~\ref{ThmUniqueLinEq} describes restrictions in potential undiscovered relations that have the form
of those in Thm.~\ref{HLaurentcoeffrel}.

In Section~\ref{SubAppnBraids} we see that when  
the first or last 
Laurent coefficient polynomial
achieves the maximum $z$-degree, this value is given by certain link invariants
in $v$ (Cors.~\ref{Cor1stProperHasMaxzdeg}, ~\ref{CorLastProperHasMaxzdeg}).
The \mbox{Homflypt} polynomials for the trivial links are unique as those for 
which the first and last
Laurent coefficient polynomials both achieve the 
maximum $z$-degree (Cor.~\ref{CorTrivialLinksUniqueness}).
When the $v$-span is 4, Thm.~\ref{Thmvspan4} displays a simple 
expression that relates the integer coefficients of  two
coefficient polynomials to the maximum $z$-degree,
the minimum $v$ degree, and the integer coefficients of a third
coefficient polynomial.
Prop.~\ref{VequivP3br} in Section~\ref{SubAppThreeBraids} establishes that 
the Jones and \mbox{Homflypt} polynomials distinguish the same three-braid links.
This depends critically on the classification of three-braids into
conjugacy classes by K.~Murasugi, \cite{70}.

The proof of Thm.~\ref{Hcoeffrel} is in an appendix as it is somewhat lengthy.

The remainder of this section has the following outline.  
Section~\ref{SubDefsAndNotation} introduces the terms and notation used in the paper.
Section~\ref{SubPropsSkeinPoly} discusses the MFW polynomials and some properties of the
\mbox{Homflypt} polynomial.
Section~\ref{SubConwayProperties}
describes properties of the Conway and Alexander polynomials.

%%%%%%%%%%%%%%%%%  Definitions and Notation   %%%%%%%%%%%%%%%%%%%%%%%%%%%%%%

\subsection{Basic Definitions and Notation}
\label{SubDefsAndNotation}
The conventions used in this paper largely follow those in,
\cite{71}, to which the reader is referred for expanded discussion.
A brief review of the standard symbols
and terminology used in the paper is given
here for reference. 

The braid group on $n$ strands, $B_n$, has $n-1$ standard generators, $\sigma_i$. 

\begin{defn}
If a braid word, $\beta \in B_n$, has the expression 
$\prod_{k=1}^m {\sigma_{i_k}^{e_k}}$, 
with $e_k = \pm 1$ for each subscript, $k$,
\begin{enumerate}%[(i)]
\item the \textit{braid length} is $m$ and is denoted $\left|\beta\right|$\,,
\item the \textit{exponent sum}, or \textit{writhe}, denoted $w$ or 
$w(\beta)$, is $\sum_{k=1}^m \,e_k$\,,
\item the mirror image of $\beta$ is $\prod_{k=1}^m {\sigma_{i_k}^{-e_k}}$, denoted
$\overline{\beta}$.
\end{enumerate}
\label{DefLengthWritheSyllable}
\end{defn}

The writhe is actually associated with the braid diagram,
rather than the braid word itself, but this imprecision in usage causes no problems.
The identity in $B_n$ is defined to have zero length and writhe.

The link associated with the 
standard closure of a braid, $\beta$, is denoted $\widehat{\beta}$.
The \textit{mirror image} of a link, $L$, is denoted $\overline{L}$. 
The number of link components in L is denoted $\mu(L)$.
$O_n$ denotes the \textit{trivial link} with $n$ components.

The skein relation, (\ref{Hskeinrel}), for the \mbox{Homflypt} polynomial, $P$, is defined
in terms of an oriented diagram, $D$, for the link.  
The diagrams, $D_+$, $D_0$, and $D_-$ below refer to the usual diagrams for the
link with positive crossing, null (smoothed) crossing, and negative crossing, respectively.

\noindent
\begin{equation}
P_{D_+}(v,z) =  vz P_{D_0}(v,z) + v^{2} P_{D_-}(v,z).
\label{Hskeinrel}
\end{equation}

The importance of the \mbox{Homflypt} polynomial derives from the fact it is 
the same for all diagrams of a given link, and so is an invariant of the oriented link.
The major skein polynomials are the \textit{Conway polynomial},
$\nabla_{D}(z) = P_{D}(1,z)$\,;
the \textit{Jones polynomial}, $V_{D}(t) = P_{D}(t, t^{\frac{1}{2}} - t^{-\frac{1}{2}})$\,;
and the \textit{Alexander polynomial}, $\Delta_{D} = P_{D}(1, t^{\frac{1}{2}} - t^{-\frac{1}{2}})$\,.
The \textit{minimum degree} of a Laurent polynomial, F, is the minimum power, 
written $\min \deg F$; for $P$ we have
$\min \deg_z P\,, \mathtt{x} = \min \deg_v P$.
The \textit{degree} of F is the highest power, 
written $\deg F$; in the
case of $P$ we have $\zeta = \deg_z P,\, \mathtt{y}=\deg_v P$.
The difference between the degree and the minimum degree is the \textit{span};
for $P$ we have $z$\textit{-span} and $v$\textit{-span}.
F is \textit{monic} when the (\textit{leading}) \textit{coefficient} for the term of maximum degree is $1$.

The \textit{torus links} are characterized by the number of full rows of twists, $p$,
together with the number of strands, $q$, denoted $K(p,q)$.
Each row may be realized by
a braid word, $\sigma_{q-1} \cdots \sigma_1$, hereafter called
$\alpha_{q-1}$, so $K(p,q)= \widehat{\alpha_{q-1}^p}$. 
The \textit{elementary torus links}, denoted $T_p$, are those with two strands. 
The corresponding Conway polynomial is so 
central to the results of this paper that the lengthy expression,
$\nabla_{T_p}(z)$, will be shortened to $C_p(z)$ or even $C_p$.

For integers $x$ we may write $\epsilon_{x}$ for $(-1)^x$. 
For $a>b$, the usual convention that 
$\sum_{j=a}^{b} f_j = 
0$, and 
$\left( \begin{array} 
{cc} 0 \\ 0
\end{array} \right) = 1 = \prod_{j=a}^{b} f_j$ 
is followed;
also 
$\left( \begin{array} 
{cc} b \\ a
\end{array} \right) = 0$ for $a <0$ or $a>b$.
The \textit{Kronecker delta}, $\delta_{a,\,b}$ is one for $a=b$ and zero otherwise.

\subsection{Properties and standard forms of the \mbox{Homflypt} polynomial}
\label{SubPropsSkeinPoly}

One way to express the \mbox{Homflypt} polynomial is as a 
polynomial in one variable, with coefficients that are 
polynomials in the other variable.
The Morton-Franks-Williams (MFW) inequality,
\cite{30}, \cite{60}, provides a 
powerful motivation to choose $v$ as the primary organizing factor, 
since it gives bounds on the possible values for the
highest and lowest powers of $v$.

Unfortunately, the MFW inequality does not tell which of
the allowed powers of $v$ have non-zero coefficients;
indeed there is no general result that describes the value of the $v$-span.
%for a given link or even when the $v$-span is sub-maximal.
However the \mbox{Homflypt} polynomial for a link
with a given braid representation ($\beta\in B_n$) may be 
organized in a standard form, (\ref{Hstdform}), in which the 
$p_j$  are ordinary (\textit{possibly zero}) polynomials with integer coefficients.
An equivalent form, (\ref{HLaurentform}), uses Laurent polynomials, $h_j = p_j/z^{n-1}$.
The $p_j$ ($h_j$) are called 
the \textit{(Laurent) MFW polynomials} for $\beta$.
When multiple braid words are under discussion, the symbols
$p_{j,\,\beta}$, or $h_{j,\,\beta}$, may be used for clarity;
note that the argument, $z$, may be omitted. 
%Both forms highlight the number of strands in the braid, and the 
%writhe, $w$, of the braid diagram (exponent sum of the braid word).
The labelling and number of MFW polynomials depends on $\beta$,
and is \textit{not} a link invariant, while  
the number of non-zero $p_{j,\,\beta}$, or $h_{j,\,\beta}$
is a link invariant.

\noindent
\begin{eqnarray}
P_{\widehat{\beta}}(v, z) &=& \frac{v^w \sum_{j=0}^{n-1} \,p_j\,(z)\,v^{2j}}
{(vz)^{n-1}}\,,   \label{Hstdform}  \\
P_{\widehat{\beta}}(v, z) &=& v^{w-n+1} \sum_{j=0}^{n-1} \,h_j(z)\,v^{2j}\,.
\label{HLaurentform}
\end{eqnarray}

\begin{defn}
Given a braid word, $\beta \in B_n$, the minimal and maximal values of j for which
$p_j \neq 0$ are designated $m(\beta)$ and $M(\beta)$ respectively; $\mathtt{m}$ and $\mathtt{M}$ may also be used
when a single braid is in use.
The (Laurent) MFW polynomials with indices in this range are called \textit{proper}.
The attributes \textit{first (last) proper} refer to polynomials of index 
$\mathtt{m}$ ($\mathtt{M}$ respectively).
The number of integers in the range $[\mathtt{m},\, \mathtt{M} ]$ is called the 
\textit{number of proper MFW polynomials},
designated $b_{\min}(\beta)$.
\label{DefMinMaxIndices}
\end{defn}
The motivation for the notation $b_{\min}(\beta)$ is that its value, $ 1 + \mathtt{M} -  \mathtt{m}$,
is well known, \cite{30}, \cite{60},
to be a lower bound for the braid index of $\widehat{\beta}$. Note that
the $v$-span of $P_{\widehat{\beta}}$ is just $\mathtt{y} - \mathtt{x} = 2(\mathtt{M} -  \mathtt{m})$.

\begin{defn}
Define $P_k(v)$ by $P_{\widehat{\beta}}(v,z) = \sum_{s \geq 0} P_{2s+1 - \mu(\widehat{\beta})}(v) z^{2s+1 - \mu(\widehat{\beta})}$.
\label{DefPk}
\end{defn}

Lemma~1, p.~133 \cite{51} (that credits Lemma~1.7 of \cite{53}) is included here:

\begin{lem}
$(v^{-1} - v)^{\mu(\widehat{\beta}) -1-s}$ divides $P_{2s+1 - \mu(\widehat{\beta})}(v)$
for $s \in [0,\,\mu(\widehat{\beta})-2]$.
\label{KawauchiLemma}
\end{lem}

The \mbox{Homflypt} formula for the elementary torus links is simply
\noindent
\begin{equation}
P_{T_p}(v, z) = \frac{v^p \{\,C_{p+1}(z) - C_{p-1}(z)\,v^{2} \}}
{(vz)^{1}}\,.
\label{Htorusstdform}
\end{equation}

%%%%%%%%%%%%%%%%%%%%%%%%%%  Properties of the Homflypt polynomial    %%%%%%%%%%%%%%%%%%%%%%%%%%
%%%%%%%%%%%%%%%%%%%%%%%%%%  Properties of the Homflypt polynomial    %%%%%%%%%%%%%%%%%%%%%%%%%%

%%%%%%%%%%%%%%%  Properties of the Homflypt polynomial for n-braid links %%%%%%%%%%%%%%%%%%%%%%%%%%%

The following result is an easy consequence of Equations \ref{Hstdform} and \ref{HLaurentform}.

\begin{prop}
The MFW polynomials obey the same skein relation (Eq.~\ref{Hskeinrel} with $v=1$) and 
reduction formulas, (\ref{Conwayrecur}), as the Conway polynomial, i.e. % $\nabla$, i.e.
\begin{eqnarray*}
h_{j,\,\gamma \sigma_i^{+1} \eta} &=& z \,h_{j,\,\gamma \sigma_i^{0} \eta} + h_{j,\,\gamma \sigma_i^{-1} \eta}\,, 
\mbox{ with } \gamma, \eta \in B_n \,, \\
%\label{Laurentskein} \\
%
h_{j,\,\beta \sigma_i^e} &=& C_e\,h_{j,\,\beta \sigma_i^{\pm 1}} + C_{e \mp 1}\,h_{j,\,\beta}\,
\mbox{ with } \beta \in B_n \,.
%\label{Laurentrecur}
\end{eqnarray*}

The MFW polynomials for $\beta$ and its mirror image, $\overline{\beta}$, satisfy
\begin{eqnarray*}
h_{j,\,  \overline{\beta}}(z) &=& h_{n-1-j,\,  \beta}(-z) = \epsilon_{w+n-1} h_{n-1-j,\, \beta}(z)\,.
\end{eqnarray*}

When $\beta \in B_n$ is extended to $\beta\sigma_{n} \in B_{n+1}$
we have $h_{n,\, \beta\sigma_{n}} = 0 = p_{n,\, \beta\sigma_{n}}$.
For all other subscripts, we have $h_{j,\, \beta\sigma_{n}} = h_{j,\, \beta}$ and
$p_{j,\, \beta\sigma_{n}} = z\,p_{j,\, \beta}$.

When $\beta \in B_n$ is extended to $\beta\sigma_{n}^{-1} \in B_{n+1}$,
we have $h_{0,\, \beta\sigma_{n}^{-1}} = 0 = p_{0,\, \beta\sigma_{n}^{-1}}$.
For all other subscripts, we have
$h_{j,\, \beta\sigma_{n}^{-1}} = h_{j-1,\, \beta}$ and
$p_{j,\, \beta\sigma_{n}^{-1}} = z\,p_{j-1,\, \beta}$.
\label{MFWcoeffprops}
\end{prop}

%%%%%%%%%%%%%%%%%%%%%%%%%%  Properties of the Conway polynomial    %%%%%%%%%%%%%%%%%%%%%%%%%%

\subsection{Properties of the Conway and Alexander polynomials}
\label{SubConwayProperties}
The Conway polynomial, $C_p$\,, is central to the results of this paper, and its 
properties are numerous and remarkable. 
The relations, $C_{p+1} = z\,C_p + C_{p-1}$, and $C_{p} = \epsilon_{p-1}\,C_{-p}$
are basic, but useful in most proofs.
As $\Delta_L(t)=\nabla_L(z=t^{\frac{1}{2}}-t^{-\frac{1}{2}})$, it is clear that any results
for $\nabla_L$ have a parallel statement for $\Delta_L$, however these will  
only be given when needed in the paper.

%%%%%%%%%%%%%%%%%% ELEM TORUS LINK CONWAY   %%%%%%%%%%%%%%%%%%%%%%%%%%%%%%%%%%%

The Conway polynomial for the elementary torus links has an expression:

\noindent
\begin{equation}
C_{p}(z) = \sum_{j=0}^{\lfloor \frac {p-1}{2} \rfloor}
\left( \begin{array} 
{cc} p-1-j \\ j
\end{array} \right) 
z^{p-2j-1}\,,\, \mbox{for } p > 0\,.
\label{Ctorusstdform}
\end{equation}

Eq.~\ref{Ctorusstdform} follows by induction and the relation of Pascal's triangle:

\noindent
\begin{eqnarray}
\left( \begin{array} 
{cc} p \\ j
\end{array} \right) &=&
\left( \begin{array} 
{cc} p-1 \\ j
\end{array} \right)
+
\left( \begin{array} 
{cc} p-1 \\ j-1
\end{array} \right) \,.
\label{Pascaltriangle}
\end{eqnarray}

Paired with $C_p$ is the Alexander polynomial. With $A_p = (t^{p} + \epsilon_{p-1})/(t+1)$

\noindent
\begin{eqnarray}
\Delta_{T_p}(t) &=& A_p/t^{(p-1)/2}\,.
\label{Deltatorusstdform} 
\end{eqnarray}
Note that $A_p = \sum_{j=0}^{p-1} \epsilon_j t^{p-1-j}$\, for $p\geq 0$, 
and $A_{-p} = \epsilon_{p-1} A_{p}/t^p$.

Use induction on $m$ and $z^{m+1}\,C_p = z^m\,C_{p+1}-z^m\,C_{p-1}$ and (\ref{Pascaltriangle}) to establish
\noindent
\begin{eqnarray}
z^m\,C_p &=& 
\sum_{j=0}^{m}
\left( \begin{array} 
{cc} m \\ j
\end{array} \right) 
\epsilon_{j}\,C_{m+p-2j}\,, \mbox{ for } m \geq 0\,.
\label{Conwayzpowers}
\end{eqnarray}
In particular, when $p=1$, there is an expression for  $z^m$.

Eq.~\ref{Conwaysumrel}-\ref{Conwaydiffofprodlast} are critical tools to 
prove the main results.
Parameters $x$, $y$, $p$, $q$, and $\kappa$, are integers with
$x+y= p+q$
in (\ref{Conwaydiffofprodlast}). Eq.~\ref{Conwaydiffofprodlast} is generally applied
with $\kappa = q$. Use induction on $y, k, \kappa$ respectively to establish (\ref{Conwaysumrel}-\ref{Conwaydiffofprodlast}).
To prove (\ref{Conwaydiffofprodlast}), also use (\ref{Conwaysumrel}) at the induction step to see 
$0=C_{x+y-2\kappa \pm 1} - C_{p+q-2\kappa \pm 1} = 
C_{x-\kappa}C_{y-\kappa} -  C_{p-\kappa}C_{p-\kappa} + C_{x-\kappa \pm 1}C_{y-\kappa \pm 1}- C_{p-\kappa \pm 1}C_{p-\kappa \pm 1}$.
\begin{eqnarray}
C_{x+y} &=& C_x\,C_{y+1} + C_{x-1}\,C_y \,, \label{Conwaysumrel}  \\
C_x\,C_{k} &=& \sum_{j=0}^{k-1} \,\epsilon_j\,C_{x+k-1-2j}\,,\mbox{ for $k \geq 0$}\,, 
\label{Conwayprodrel}  \\
C_x\,C_y -C_p\,C_q     &=& 
\epsilon_{\kappa}\, \{ C_{x-\kappa}\,C_{y-\kappa} -
C_{p-\kappa}\,C_{q-\kappa}
\}, \mbox{ for } x+y= p+q\,.
\label{Conwaydiffofprodlast}
\end{eqnarray}

\begin{prop}
There are no common roots over the complex numbers for $C_p$ and
$C_{p+1}$, for any integer $p \neq 0, -1$, hence $\gcd (C_p\,,\,C_{p+1}) = 1$\,.

Second, for any integers $a \neq 0$, and $b$, we have $C_a |\, C_{ab}$\,.

Third, when $\gcd (a,b) = g$, we have $\gcd (C_a,C_b) = C_g$\,.
\label{Conwayidealprop}
\end{prop}
\begin{proof}
It suffices to consider only positive integers $p, a, b$ as $C_0=0$ and
$C_x = \epsilon_{x+1} C_{-x}$ for $x < 0$.
The first claim follows from the relation, $C_{p+1}=z\,C_p + C_{p-1}$, 
as any common root of $C_p$ and $C_{p+1}$ would then be a root of of $C_{p-1}$ and
by induction a root of $C_1=1$.
For the second claim, observe that $C_{ab}=C_a C_{ab-a+1} +
C_{a-1} C_{ab-a}$, by (\ref{Conwaysumrel}), and use induction on $b$.
For the third claim, there are integers, $\lambda$, $\mu$ so that
$\lambda a + \mu b = g$. Eq.~\ref{Conwaysumrel} now implies $C_g = C_{\lambda a + \mu b} = 
C_{\lambda a + 1}C_{\mu b} + C_{\lambda a}C_{\mu b-1}$; by 
the second claim any common root of $C_a$ and $C_b$ is a root of $C_g$. 
The second claim also gives us $C_g |\, C_a, C_b$\,.
\end{proof}

%%%%%%%%%%%%

Use induction on $e$ to prove
(\ref{Conwayrecur}) that shows how the Conway polynomial can be 
calculated using braid words of shorter length.
For $\beta \in B_n$  and any integer, $e$, we have
\noindent
\begin{equation}
\nabla_{\widehat{\beta \sigma_i^e}} = 
C_e\,\nabla_{\widehat{\beta \sigma_i^{\pm 1}}} +
C_{e \mp 1}\,\nabla_{\widehat{\beta}}\,.
\label{Conwayrecur}
\end{equation}
%%%%%%%%%%%%%%%  Properties of the Conway polynomial for three-braid links %%%%%%%%%%%%%%%%%%%%%%%%%%%%%
The following proposition is useful in the applications to three-braid link results in the paper.
Eq.~\ref{Alexander3brApprox} is essentially given by Prop.~4.2, p.~13-14 \cite{70}.
\begin{prop}
%When $\gamma \in B_3$ and $a,\,x$ are integers with $a \geq 0$, we have:
When $\gamma \in B_3$ and $a,\,x \in \mathbb{Z}$ with $a \geq 0$ and $G = t^2+t+1$:
\begin{eqnarray}
\nabla_{\widehat {\alpha_2^3 \gamma} } &=& \nabla_{\widehat {\gamma} } + C_{w(\gamma)+5} - C_{w(\gamma)+1} \,, 
\label{ConwayBasicCenter}  \\
\nabla_{\widehat {\alpha_2^{3a} \gamma} } &=& 
\nabla_{\widehat {\gamma} }+ \sum_{j=1}^a C_{w(\gamma)+6j-1} - \sum_{j=1}^a C_{w(\gamma)+6j-5}\,,  
\label{ConwayPosPowerCenter}  \\
\nabla_{\widehat {\alpha_2^{-3a} \gamma} } &=& 
\nabla_{\widehat {\gamma} }+ \sum_{j=1}^a C_{w(\gamma)-6j+1} - \sum_{j=1}^a C_{w(\gamma)-6j+5}\,,
\label{ConwayNegPowerCenter} \\
\Delta_{\widehat {\alpha_2^{3x} \gamma} } &=&
\Delta_{\widehat {\gamma} }+ t ( t^{w(\gamma)+3x}  - \epsilon_{w(\gamma)} )(t^{3x}-1)/(Gt^{3x}t^{w(\gamma)/2})\,. \label{DeltaAnyPowerCenter} 
\end{eqnarray}

For $\gamma_r = \prod_{k=1}^{r} \sigma_2^{-e_{k,2}} \sigma_1^{e_{k,1}}$, with $r, e_{k,2}\,, e_{k,1} >0$, and
$1 < E_i=\sum_{k=1}^r e_{k,i}$, 
\begin{eqnarray}
\nabla_{\widehat {\gamma_r} } &=& (-1)^{E_2+1} \{C_{E_2+E_1-1} - r C_{E_2+E_1-3} + o(C_{E_2+E_1-3})    \}\,,
\label{Conway3brApprox} \\
\Delta_{\widehat {\gamma_r} } &=& (-1)^{E_2+1} \{t^{(E_2+E_1-2)/2} - (r+1) t^{(E_2+E_1-4)/2} + \cdots    \}\,.
\label{Alexander3brApprox}
%
%\nabla_{\widehat {\gamma} } &=& (-1)^{E_2+1} \{C_{E_2} C_{E_1} - \prod_{k=1}^{2} C_{e_{k,2}} C_{e_{k,1}} \}\,, \mbox{ when $r=2$.} \nonumber 
\end{eqnarray}

\label{Conway3braugmented}
\end{prop}

\begin{proof}
Eq.~\ref{ConwayBasicCenter}: use $\alpha_2^3 = (\sigma_1^2 \sigma_2)^2$ and (\ref{Conwayrecur}) to see
$\nabla_{\widehat {\alpha_2^3 \gamma} }= z \nabla_{\widehat { \sigma_1^2 \sigma_2 \sigma_1 \sigma_2\gamma} }
+ z \nabla_{\widehat { \sigma_1 \sigma_2^2 \gamma} } + z \nabla_{\widehat { \sigma_2\gamma} } + \nabla_{\widehat {\gamma} }$.
Since $C_{w(\gamma)+5} - C_{w(\gamma)+1} = zC_{w(\gamma)+4} + z C_{w(\gamma)+2}$, it is equivalent to show
$ C_{w(\gamma)+4} +  C_{w(\gamma)+2} = \nabla_{\widehat { \sigma_1^2 \sigma_2 \sigma_1 \sigma_2\gamma} }
+ \nabla_{\widehat { \sigma_1 \sigma_2^2 \gamma} } +  \nabla_{\widehat { \sigma_2^1\gamma} }$, and this is 
the same as 
$C_{w(\gamma)+4} +  C_{w(\gamma)+2} =
\nabla_{\widehat { \sigma_2 \sigma_1 \sigma_2^3 \gamma} }
+ \nabla_{\widehat { \sigma_1 \sigma_2^2 \gamma} } +  \nabla_{\widehat { \sigma_2^1\gamma} }$\,.

When $\gamma = \sigma_2^e \sigma_1^f$, note that 
$\widehat { \sigma_2 \sigma_1 \sigma_2^3 \gamma} = \widehat { \sigma_2 \sigma_1 \sigma_2^{e+3} \sigma_1^f}
= \widehat { \sigma_1^{e+3} \sigma_2 \sigma_1 \sigma_1^f} = T_{e+f+4}$.
Also $\nabla_{\widehat { \sigma_1 \sigma_2^2 \gamma} } +  \nabla_{\widehat { \sigma_2\gamma} }=
C_{e+2} C_{f+1} + C_{e+1} C_{f} $, and this is just $C_{e+f+2}$ by (\ref{Conwaysumrel}). 
This shows (\ref{ConwayBasicCenter})
holds for $\gamma = \sigma_2^e \sigma_1^f$.

Assume now that (\ref{ConwayBasicCenter}) is true for $\gamma_r = \prod_{i=1}^{r} \sigma_2^{e_i} \sigma_1^{f_i}$.
The following observation facilitates the induction step:
$\widehat {  \gamma_r \sigma_2 \sigma_1} = \widehat {  \gamma_r^{*} }$, with 
$e_1^{*} = e_1+f_1+1$ and $f_1^{*} = 1$ when $r=1$; 
for $r>1$, use
$e_1^{*} = e_1+f_r$ and $e_r^{*}= e_r+1$ and $f_r^{*} = 1$, and otherwise set $e_j^{*} = e_j$ for $1<j<r$ and 
$f_j^{*} = f_j$ for $j \leq 1<r$.

Let $R=r+1$ and use (\ref{Conwayrecur}) to see that 
$\nabla_{\widehat {\alpha_2^3 \gamma_R} }$ is just the sum of four terms:
$C_{e_R} C_{f_R} \nabla_{\widehat {\alpha_2^3 \gamma_r \sigma_2 \sigma_1} }$;
$C_{e_R} C_{f_R-1} \nabla_{\widehat {\alpha_2^3 \gamma_r \sigma_2 } }$;
$C_{e_R-1} C_{f_R} \nabla_{\widehat {\alpha_2^3 \gamma_r \sigma_1 } }$;
$C_{e_R-1} C_{f_R-1} \nabla_{\widehat {\alpha_2^3 \gamma_r } }$.
A simple manipulation of the braids involved allows us to apply the
induction hypothesis to obtain the sum of: 
$C_{e_R} C_{f_R} (\nabla_{\widehat { \gamma_r \sigma_2 \sigma_1} } + C_{w(\gamma_r)+7} -  C_{w(\gamma_r)+3}) $; 
$C_{e_R} C_{f_R-1} (\nabla_{\widehat {\gamma_r \sigma_2 } } + C_{w(\gamma_r)+6} -  C_{w(\gamma_r)+2}) $; 
$C_{e_R-1} C_{f_R} (\nabla_{\widehat {\gamma_r \sigma_1 } } + C_{w(\gamma_r)+6} -  C_{w(\gamma_r)+2}) $;
$C_{e_R-1} C_{f_R-1} (\nabla_{\widehat {\gamma_r } } + C_{w(\gamma_r)+5} -  C_{w(\gamma_r)+1}) $.
Use Eqs.~\ref{Conwayrecur}, \ref{Conwaysumrel} to combine the first pair to obtain
$C_{e_R} (\nabla_{\widehat { \gamma_r \sigma_2 \sigma_1^{f_R}} } + C_{w(\gamma_r)+f_R+ 6} -  C_{w(\gamma_r)+f_R+2}) $
and then combine the second pair to obtain
$C_{e_R-1} (\nabla_{\widehat {\gamma_r \sigma_1^{f_R} } } + C_{w(\gamma_r)+f_R+ 5} +  C_{w(\gamma_r)+f_R+1}) $.
A final use of (\ref{Conwayrecur}, \ref{Conwaysumrel}) gives us 
$\nabla_{\widehat { \gamma_r \sigma_2^{e_R} \sigma_1^{f_R}} } + C_{w(\gamma_r)+e_R+f_R+ 5} -  C_{w(\gamma_r)+e_R+f_R+1} $, as desired.

Eqs.~\ref{ConwayPosPowerCenter} and \ref{ConwayNegPowerCenter} follow from (\ref{ConwayBasicCenter}) by induction.
To verify Eq.~\ref{DeltaAnyPowerCenter} for $x >0$, use (\ref{ConwayPosPowerCenter}), convert
$\nabla_L$ to $\Delta_L$, and use (\ref{Deltatorusstdform}) 
with $\Delta_{T_p}= (t^{(1+p)/2} + \epsilon_{p-1} t^{(1-p)/2} )/(t+1)$
to see each summation in (\ref{ConwayPosPowerCenter}) as two geometric series and combine.
When $x < 0$, use $\Delta_{L}(t) = \Delta_{\overline{L}}(t^{-1})$,
so $\Delta_{\widehat {\alpha_2^{3x} \gamma} }(t) = \Delta_{\widehat {\alpha_2^{3|x|} \overline{\gamma}}} (t^{-1}) $.

Eq.~\ref{Alexander3brApprox} is essentially Prop.~4.2, p.~13-14 \cite{69}, and equivalent to
Eq.~\ref{Conway3brApprox}. Use (\ref{Conwayprodrel}) to prove (\ref{Conway3brApprox}) for $r=1$.
For $r=2$ a proof by induction on the length of $\gamma_2$ shows
$\nabla_{\widehat {\gamma_2} } = (-1)^{E_2+1} \{C_{E_2} C_{E_1} - \prod_{k=1}^{2} C_{e_{k,2}} C_{e_{k,1}} \}$
that satisfies (\ref{Conway3brApprox}).
Set $R=r+1$ for $r \geq 2$, and apply (\ref{Conwayrecur}) at $\sigma_2^{-e_{R,2}}$ and $\sigma_1^{e_{R,1}}$
so as to express $\nabla_{\widehat {\gamma_R} }$ as a sum of four terms associated with
alternating braids. 
Apply the skein relation to $\gamma_r \sigma_2^{-1} \sigma_1^1$ to replace
$C_{-e_{R,2}} C_{e_{R,1}} \nabla_{\widehat{\gamma_r \sigma_2^{-1} \sigma_1^1} }$
by two new terms related to $\gamma_r \sigma_2^{+1} \sigma_1^1$ 
and $\gamma_r \sigma_2^{0} \sigma_1^1$.
As above, replace the first new term for 
$\gamma_{r-1}  \sigma_2^{-e_{r,2}}\sigma_1^{e_{r,1}}  \sigma_2^{+1}\sigma_1^{1}$, by 
a term for $\gamma_{r-1}  \sigma_2^{-e_{r,2}}\sigma_2^{1}  \sigma_1^{+1}\sigma_2^{e_{r,1}}$, i.e.
$\xi=\sigma_2^{e_{r,1}}\gamma_{r-1}  \sigma_2^{1-e_{r,2}}  \sigma_1^{+1}$.
Set $d_k = - e_{k,2}$ and $f_k= e_{k,1}$.
Use (\ref{Conwayrecur}) to replace the $\xi$ term by four terms related to alternating braids;
these appear as the first two lines in (\ref{Conway3brApproxproof1}).
Rewrite $C_{d_R+1} C_{f_R} \nabla_{\widehat {\gamma_{r} \sigma_1 } } - z C_{d_R } C_{f_R} \nabla_{\widehat {\gamma_{r} \sigma_1 } }$
% as $ C_{d_R-1} C_{f_R}  \nabla_{\widehat {\gamma_{r} \sigma_1 } }$ 
to see $\nabla_{\gamma_{R}}$ can now be expressed 
as a sum of terms to which the
induction hypothesis applies:
\begin{eqnarray}
\nabla_{\gamma_{R}}
&=&  C_{d_R} C_{f_R} C_{f_r} z \nabla_{\widehat {\sigma_2^{-1} \gamma_{r-1}  \sigma_2^{d_r}   \sigma_1 } } % 1st
+  C_{d_R} C_{f_R} C_{f_r}  \nabla_{\widehat {\sigma_2^{-1} \gamma_{r-1}  \sigma_2^{d_r-1}  \sigma_1 } } \nonumber \\ % 2nd
&&+  C_{d_R} C_{f_R} C_{f_r+1} z \nabla_{\widehat {\gamma_{r-1}  \sigma_2^{d_r}   \sigma_1  } }% 3rd
+ C_{d_R} C_{f_R}  C_{f_r+1} \nabla_{\widehat {\gamma_{r-1}  \sigma_2^{d_r-1}  \sigma_1  } } \nonumber \\ %4th
&&+  C_{d_R-1} C_{f_R}  \nabla_{\widehat {\gamma_{r} \sigma_1 } }% 5th
+  C_{d_R} C_{f_R-1}  \nabla_{\widehat {\sigma_2^{-1} \gamma_{r} } }% 6th
+  C_{d_R+1} C_{f_R-1}  \nabla_{\widehat {\gamma_{r} } }\,.  % 7th
\label{Conway3brApproxproof1}
\end{eqnarray}
Apply the induction hypothesis to each term and
use Eqs.~\ref{Conwaysumrel}, \ref{Conwaydiffofprodlast} and 
the skein relation to simplify and obtain (\ref{Conway3brApprox}) for
$\nabla_{\gamma_{R}}$.
\end{proof}

%%%%%%%%%%%%%%%%%%%%%%%%%%%%%%%%%%%%%%%%%%%%%%%%%%%%%%%%%%%%%%%%%%%%%%%%%%%

\section{The Primary Results}
\label{SectionMain}

In 1987, V.F.R.~Jones published a number of results relating properties
of Hecke algebras to knot theory, \cite{41}. This landmark research displayed, among
other insights, how the two variable skein polynomial could be calculated
from the Burau representation of a braid. 
%Unfortunately for knot theorists, the expression for
%the two variable skein polynomial is of nontrivial complexity.  
For three-braid links, this paper first established the
\mbox{Homflypt} polynomial is dependent only on the writhe 
and the Alexander polynomial.

H.~Murakami, \cite{69}, shows how to calculate $P_{D}(a,z)$
for a diagram with $n$ Seifert circles based on knowledge of $P_{D}(a_j,z)$ 
at $n$ independent points $a_j$. A corollary, \cite{69}, shows that
for diagrams with 3, 4, or 5 Siefert circles
$P_L(a,t^{\frac{1}{2}} - t^{-\frac{1}{2}})$ may be expressed in terms of weighted sums of $1$,
$\Delta_L(t)$, $V_L(t)$, and $V_L(t^{-1})$, where the weights are certain
Laurent polynomials
in $a, \sqrt t $.
However, there is no expression for $P_L(a,z)$ in terms of $a, z$ and directly
inverting the expressions in $a, \sqrt t$ to become expressions in $a, z$ appears
essentially impossible for the general case.

Theorem~\ref{Hcoeffrel} provides an expression for $P_L(v,z)$ in terms of $v, z$.
The first three MFW 
polynomials, $h_0(z)$, $h_1(z)$, and $h_2(z)$, which may be zero, are determined by the writhe,
Conway polynomial, and the $h_j(z)$ for $j \geq 3$. 
%Recall that some of these MFW polynomials may be zero. 
A similar expression for the proper MFW polynomials 
appears in Cor.~\ref{HProperLcoeffrel}.
 
\begin{thm}
For an arbitrary braid, $\beta$, of $n \geq 1$ strands, the \mbox{Homflypt} polynomial for 
$\widehat{\beta}$ is given by (\ref{HLaurentform}) with:
\begin{enumerate}%[(i)]
\item $h_0 = z^{-2}\{\,C_{w+4-n} - \nabla_{\widehat{\beta}}\} -  q_{0} $\,,
with $q_{0} = \sum_{j=3}^{n-1}  \,z^{-2}(C_{2j-3}-1) \,h_j$\,,
\item $h_{1} = \nabla_{\widehat{\beta}} -h_0 - h_2 - q_{1} $\,,
with $q_{1} = \sum_{j=3}^{n-1} \,h_j$\,,
\item $h_2 = z^{-2}\{\,C_{w+2-n} - \nabla_{\widehat{\beta}}\}  - q_{2}  $\,,
with $q_{2} = \sum_{j=3}^{n-1}  \,\,z^{-2}(C_{2j-1}-1) \,h_j$\,.
\end{enumerate}

\label{Hcoeffrel}
\end{thm}

A few comments are in order about this theorem. First, the formula for 
$n=1$ depends on a writhe of zero. 
Second, the formula for $n=2$ is equivalent to (\ref{Htorusstdform}).
Third, the formula for $n \leq 3$ has $0 = q_0 = q_1 = q_2$.
Finally, $h_{1}$ is exactly what is needed to satisfy the identity, 
$P_{\widehat{\beta}}(1, z) = \nabla_{\widehat{\beta}}(z)$;
$h_1$ could be rewritten in a fashion similar to $h_0, h_2$.
This result appears to be a consequence of H.~Murakami's Theorem 1, p.~409 \cite{69}, but 
the suitable choice of $n$ independent values is unclear; neither will the Corollary,
p.~410 \cite{69}, suffice for the $n$-braid case.
See Appendix~\ref{ProofHcoeffrel} for a proof.

A closer inspection of the forms for $h_j$, $q_0$, and 
$q_2$, reveals a set of relations, described by  
Theorem~\ref{HLaurentcoeffrel} (recall some $h_j$ may be zero).
Independent means that for a fixed choice of 
$w$ and $n$, a set of $ f_j $  can satisfy one of (\ref{HLaurentrel1}, \ref{HLaurentrel2})
and not satisfy the other. 

%%%%%%%%%%%%%%%%%%%%%%%% LINEAR EQ Homflypt %%%%%%%%%%%%%%%%%%%%%%%%%%%%%%%%%%%%%%%%%

\begin{thm}
The Laurent MFW polynomials associated with a braid, $\beta \in B_n$,
satisfy the following relations with $f_j=h_j$. 
These relations are independent for $n>1$\,.
 
\noindent
\begin{eqnarray}
\sum_{j=0}^{n-1}   \,C_{2j-3}\,f_j &=&  C_{4+w-n} \label{HLaurentrel1}\,, \\
\sum_{j=0}^{n-1}   \,C_{2j-1}\,f_j &=&  C_{2+w-n} \label{HLaurentrel2}\,.
\end{eqnarray}

Any set of $n$ functions, $\{f_j\}_0^{n-1}$, 
that satisfy (\ref{HLaurentrel1}) and (\ref{HLaurentrel2}), will also satisfy the following family of equations
for each integer, $\kappa$:
\noindent
\begin{equation}
\sum_{j=0}^{n-1}   \,C_{2j-1-\kappa }\,f_j =  \epsilon_{\kappa}\,C_{\kappa+2+w-n}\,.
\label{HLaurentrel3}
\end{equation}

\label{HLaurentcoeffrel}
\end{thm}

\begin{proof}
To prove (\ref{HLaurentrel1}, \ref{HLaurentrel2}) substitute 
the values of $h_j$ from Thm.~\ref{Hcoeffrel} for $f_j$.

To show (\ref{HLaurentrel1}), (\ref{HLaurentrel2})
are independent when $n>1$, it suffices to find two sets of
$n$ polynomials each of which satisfies one relation, but not the other.

In case $w\neq n-1$, Eq.~\ref{Conwaysumrel} implies $f_0=C_{w+2-n}$, $f_1=z\,C_{w+1-n}$, and
$f_j=0$, for $j>1$ satisfy (\ref{HLaurentrel1}), but these 
don't satisfy (\ref{HLaurentrel2}). When $w=n-1$, the choice $f_1=C_3$ and $f_j=0$, for
$j\neq 1$ satisfies (\ref{HLaurentrel1}), but not (\ref{HLaurentrel2}).

In case $w\neq n-3$, observe that $f_1=C_{w+2-n}$, and 
$f_j=0$, for $j\neq 1$ satisfy (\ref{HLaurentrel2}), but they
don't satisfy (\ref{HLaurentrel1}). When $w=n-3$, the choice $f_0=1$, and $f_j=0$, for
$j\neq 0$ satisfies (\ref{HLaurentrel2}), but does not satisfy (\ref{HLaurentrel1}).

%%%%%%%%%%%%%%%%%%%

The final point is to show whenever a function satisfies 
(\ref{HLaurentrel1}, \ref{HLaurentrel2}), 
it also satisfies (\ref{HLaurentrel3}). First take the difference of 
(\ref{HLaurentrel1}, \ref{HLaurentrel2}), and divide by $z$ to obtain (\ref{HLaurentrel3}) with $\kappa = 1$. 
Now multiply both sides of Eq.~\ref{HLaurentrel3} with $\kappa = 1$, by $C_{\kappa}$,
multiply both sides of (\ref{HLaurentrel2}) by $-C_{\kappa-1}$, and add these 
two products. 
Apply (\ref{Conwaysumrel}) to the right side of this sum and 
(\ref{Conwaydiffofprodlast}) to the left side to derive (\ref{HLaurentrel3}).
\end{proof}

Theorem~\ref{HLaurentcoeffrel} leads to an expression for the 
first three proper Laurent MFW polynomials in terms of
the remaining $h_j$ and three
link invariants: $\mathtt{x} = \min \deg_v P_{\widehat{\beta}}$\,, 
$v$-span $P_{\widehat{\beta}}$\,, and $\nabla_{\widehat{\beta}}$\,.

\begin{cor}
Given an arbitrary braid, $\beta$, of $n \geq 1$ strands, 
for which
$\mathtt{M} \geq \mathtt{m}+2$, the \mbox{Homflypt} polynomial for 
$\widehat{\beta}$ is given by (\ref{HLaurentform}) with:
\begin{enumerate}%[(i)]
\item $h_{\mathtt{m}} = z^{-2} \{ C_{\mathtt{x}  +3} -  \nabla_{\widehat{\beta}}  
-  \sum_{j=3}^{\mathtt{M}-\mathtt{m}} ( C_{2j-3 } -1) h_{\mathtt{m}+j} \} $\,,
\label{HProperLcoeffrel0}
\item $h_{\mathtt{m}+1} = \nabla_{\widehat{\beta}} - h_{ \mathtt{m} } - h_{\mathtt{m}+2} 
- \sum_{j=3}^{\mathtt{M}-\mathtt{m}} \,h_{\mathtt{m}+j}$\,, which is true by the definitions,
\label{HProperLcoeffrel1}
\item $h_{\mathtt{m}+2}  = z^{-2} \{ C_{\mathtt{x}  +1} - \nabla_{\widehat{\beta}}  
-\sum_{j=3}^{\mathtt{M}-\mathtt{m}} (C_{2j-1 } - 1) h_{\mathtt{m}+j} \}$\,.
\label{HProperLcoeffrel2}
\end{enumerate}

For $\mathtt{M} < \mathtt{m}+2$, replace $h_{\mathtt{m}+2}$ by $0$ in (\ref{HProperLcoeffrel2}).
For $\mathtt{M}=\mathtt{m}$ we have $P_{\widehat{\beta}}=1$. 
For $\mathtt{M}=\mathtt{m}+1$, we have $P_{\widehat{\beta}}=v^{\mathtt{x}  } (C_{\mathtt{x}  +2}- v^2 C_{ \mathtt{x} }  )/z$. 
\label{HProperLcoeffrel}
\end{cor}

\begin{proof}
We see $h_{\mathtt{m}} + h_{\mathtt{m}+1} + C_3 h_{\mathtt{m}+2} 
+ \sum_{j=\mathtt{m}+3}^{n-1} C_{2j-2\mathtt{m}-1 }h_{j} =  C_{\mathtt{x}  +1}$ by use of 
$\kappa = 2 \mathtt{m}$ in (\ref{HLaurentrel3}).
Now
$z^2 h_{\mathtt{m}+2}  =  C_{\mathtt{x}  +1} - \nabla_{\widehat{\beta}}  -\sum_{j=\mathtt{m}+3}^{n-1} [C_{2j-2\mathtt{m}-1 } - 1]h_{j}$
follows by substitution for $h_{\mathtt{m}+1}$ from (\ref{HProperLcoeffrel1}).
This establishes (\ref{HProperLcoeffrel2}).

Also we have $h_{\mathtt{m}} + C_3 h_{\mathtt{m}+1} + C_5 h_{\mathtt{m}+2} 
+ \sum_{j=\mathtt{m}+3}^{n-1} C_{2j-2\mathtt{m} +1 }h_{j} =  C_{\mathtt{x}  -1}$
by use of $\kappa = 2\mathtt{m}-2$ in (\ref{HLaurentrel3}). 
Now use the value of $h_{\mathtt{m}+1}$ from (\ref{HProperLcoeffrel1}) to obtain
$$z^2 h_{\mathtt{m}} = C_3 \nabla_{\widehat{\beta}} 
+ (C_5 - C_3) h_{\mathtt{m}+2} + \sum_{j=\mathtt{m}+3}^{n-1} [C_{2j-2\mathtt{m} +1 } - C_3  ] h_{j} -  C_{\mathtt{x}  -1}\,.$$
Observe that $C_5 - C_3$ is just $z^2(z^2+2)=z^2(C_3+1)$.
Use (\ref{HProperLcoeffrel2}) 
to express
$z^2 h_{\mathtt{m}+2}$ and simplify to obtain 
\begin{eqnarray*}
z^2 h_{\mathtt{m}} &=& (C_3+1) C_{\mathtt{x}  +1} -  C_{\mathtt{x}  -1} -  \nabla_{\widehat{\beta}}  \\ 
&&-  \sum_{j=\mathtt{m}+3}^{n-1} (C_3+1)[C_{2j-2\mathtt{m}-1 } - 1]h_{j}
+ \sum_{j=\mathtt{m}+3}^{n-1} [C_{2j-2\mathtt{m} +1 } - C_3  ] h_{j}\,.
\end{eqnarray*}

Use (\ref{Conwaysumrel}) to see that $C_{\mathtt{x}  +3}= C_3 C_{\mathtt{x}  +1} + z C_{\mathtt{x} } $, 
that is $(C_3+1) C_{\mathtt{x}  +1} -  C_{\mathtt{x}  -1}$.
Similarly, the coefficient of $h_{j}$ is just
$1-C_{2j-2\mathtt{m}-1 } +z  C_{2j-2\mathtt{m} -2 }$, i.e. $1-C_{2j-2\mathtt{m}-3 }$.
Putting this together we have the following relation equivalent to (\ref{HProperLcoeffrel0}).
\begin{eqnarray*}
z^2 h_{\mathtt{m}} &=& C_{\mathtt{x}  +3} -  \nabla_{\widehat{\beta}}  -  \sum_{j=\mathtt{m}+3}^{n-1} ( C_{2j-2\mathtt{m}-3 } -1) h_{j}\,.
\end{eqnarray*}
%$\sum_{j=0}^{n-1}   \,C_{2j-1-\kappa }\,f_j =  \epsilon_{\kappa}\,C_{\kappa+2+w-n}\,.$

When $\mathtt{M}=\mathtt{m}$, we have  $C_{\mathtt{x}  +1} = \nabla_{\widehat{\beta}}$ (by \ref{HProperLcoeffrel2}),
so (\ref{HProperLcoeffrel1}) tells us $h_{ \mathtt{m} } = C_{\mathtt{x}  +1}$, and (\ref{HProperLcoeffrel0}) now implies
$h_{ \mathtt{m} } = z^{-2} (C_{\mathtt{x}  +3} - C_{\mathtt{x}  +1}) = C_{\mathtt{x}  +1}$. 
This is equivalent to $z^{-1} C_{\mathtt{x}  +2} =  C_{\mathtt{x}  +1}$ whose only solution is $\mathtt{x}  =0$.
Thus $h_{ \mathtt{m} }=1 =P_{\widehat{\beta}}$\,. The case $\mathtt{M}=\mathtt{m}+1$ is straightforward to prove.
\end{proof}

Given the number of possible applications (Section~\ref{SubAppnBraids}), 
it would clearly be desirable to find other independent relations.
Suppose some equation exists in the form of Theorem~\ref{HLaurentcoeffrel}
for each value of $n$, with coefficients, $\theta_j$, that are independent of $\beta$ and $n$.
We want $\Omega_{n,\,\beta}$ to be expressed in terms of the
braid properties of $\beta$, as $\Omega_{n,\,\beta}$ could always be trivially \textit{defined}
by Eq.~\ref{Hconjformula}.

\noindent
\begin{eqnarray}
\sum_{j=0}^{n-1}   \,\theta_{j}\,h_j &=&  \Omega_{n,\,\beta}\,, \mbox{ with }  \theta_{j} \in \mathbb{Z}{[z, z^{-1}]}\,.
\label{Hconjformula}
\end{eqnarray}

We already know that 
$\theta_j=1$ works with $\Omega_{n,\,\beta} = \nabla_{\widehat{\beta}}$, as 
does $\theta_j=C_{2j-3}$ with $\Omega_{n,\,\beta} = C_{4+w-n}$, or 
$\theta_j=C_{2j-1}$ with $\Omega_{n,\,\beta} = C_{2+w-n}$.
Observe that $\nabla_{\widehat{\beta}}$ is a link invariant, while
$C_{4+w-n}$ and $C_{2+w-n}$ are not.
On the other hand, $C_{4+w-n}$ and $C_{2+w-n}$ are functions of $w-n$,
while $\nabla_{\widehat{\beta}}$ is not.
Furthermore, for any $a_0, a_1, a_2 \in \mathbb{Z}{[z, z^{-1}]}$, the weighted sum,
$\Omega_{n,\,\beta}^{(a_0, a_1, a_2)} = a_0 C_{4+w-n} + a_1 \nabla_{\widehat{\beta}}+ a_2 C_{2+w-n}$ satisfies (\ref{Hconjformula})
with $\theta_j = a_0 C_{2j-3} + a_1 + a_2 C_{2j-1}$.
When $a_0, a_1, a_2$ are all non-zero, $\Omega_{n,\,\beta}^{(a_0, a_1, a_2)}$ is neither a link invariant, 
nor a function of $w-n$.
Are there other, independent, choices for $\Omega_{n,\,\beta}$?

As any two conjugate braid words in $B_n$ produce
the same set of $h_j$, the value of $\Omega_{n,\,\beta}$ must be identical for them.
Any pair of braid words that 
generate the same braid element in $B_n$ also have this same property,  so
$\Omega_{n,\,\beta}$ is a function on the conjugacy classes of braid elements. 
In fact, by (\ref{HLaurentform}), any two braid elements, even in different braid groups, 
with identical $w_i-n_i$ and whose closure is the
same link, will have $\Omega_{n_1,\,\beta_1} = \Omega_{n_2,\,\beta_2}$.
By Prop.~\ref{MFWcoeffprops},
the first $n$ MFW polynomials, $h_j$, are identical for 
$\beta \in B_n$ and $\beta \sigma_n \in B_{n+1}$, and $h_{n,\, \beta \sigma_n}=0$,
so the left side of (\ref{Hconjformula})
is \textit{invariant under positive Markov stabilization}, i.e. 
$\Omega_{n,\,\beta} =  \Omega_{n+1,\,\beta\sigma_n}$.

The trivial knot has multiple representations of the form 
$\prod_{k=1}^{n-1} {\sigma_{k}^{e_k}}$, with $e_k= \pm 1$, for a given writhe between $1-n$ and $n-1$.
Critically, all such representations
with the same writhe have the same set of MFW polynomials, $h_j$. 
Observe that any such choice of braid word, $\omega_j$, with 
$j$ negative exponents, $e_k$, and $n-1-j$ positive exponents,
satisfies $\theta_j = \Omega_{n,\,\omega_j}$\,, since $h_{s,\, \omega_j} = \delta_{j,\,s}$\,. 
If $\Omega_{n,\,\beta}$ is invariant under negative Markov stabilization,
$\Omega_{n,\,\beta}$ is a link invariant. In this case,
observe that for the trivial knot, $\theta_0 = \Omega_{2,\,\sigma_1}= \Omega_{n,\, \omega_j}=\theta_j$, 
Thus (\ref{Hconjformula}) becomes $\theta_0~\nabla_{\widehat{\beta}} = \Omega_{n,\,\beta}$,
so $\Omega_{n,\,\beta}$ is a multiple of $\nabla_{\widehat{\beta}}$.

Eq.~\ref{Htorusstdform} shows that for elementary torus links, we have
$h_{0,\, \sigma_1^{\nu+2}} = C_{\nu+3}/z$ and $h_{1,\, \sigma_1^{\nu+2}} = -C_{\nu+1}/z$, and
by (\ref{Hconjformula}):
\begin{eqnarray}
\{\theta_0\,C_{\nu+3} -\theta_1\,C_{\nu+1}\}/z = \Omega_{2,\,\sigma_1^{\nu+2}} \,.
\label{GenlTrel}
\end{eqnarray}

Suppose $\Omega_{n,\,\beta}$ depends only on $w-n$, 
as is true when $\beta=\prod_{k=1}^{n-1} {\sigma_{k}^{e_k}}$, with $e_k= \pm 1$. 
Then for any $\beta \in B_n$, we would have 
$\Omega_{n,\,\beta} = \Omega_{2,\,\sigma_1^{w-n+2}}$
and $\theta_j= \Omega_{n,\,\omega_j}= \{ \theta_0\,C_{-2j+2} - \theta_1\,C_{-2j} \}/z$ 
(note $\nu=w-n=-1-2j$ for $\omega_j$). 
Substitution of these values in (\ref{Hconjformula}) and use of $\kappa=1, -1$ in Eq.~\ref{HLaurentrel3}
reveals an identity, i.e. any $\theta_0, \theta_1 \in \mathbb{Z}{[z, z^{-1}]}$ are valid.
In other words, 
$\Omega_{n,\,\beta}= \theta_0 C_{3+w-n}/z - \theta_1 C_{1+w-n}/z$ is merely a weighted sum of functions in Thm.~\ref{HLaurentcoeffrel},
as $C_{3+w-n} = (C_{4+w-n} - C_{2+w-n})/z$ and $C_{1+w-n} = (C_{4+w-n} - C_3 C_{2+w-n})/z$.

The question that remains is whether any $\Omega_{n,\,\beta}$ that satisfies (\ref{Hconjformula})
is merely a weighted sum of the functions we have already discovered. Observe that when 
$\Omega_{n,\,\beta}^{*}$ and $\Omega_{n,\,\beta}^{**}$ satisfy (\ref{Hconjformula}), so also will
$\Omega_{n,\,\beta}^{*} \pm \Omega_{n,\,\beta}^{**}$.
Starting from an arbitrary such $\Omega_{n,\,\beta}$, we know that
$\theta_j = \Omega_{n,\,\omega_j}$. 
Now consider  $\Omega_{n,\,\beta}^{*} = a_0 C_{w+4-n} + a_1 \nabla_{\widehat{\beta}} + a_2 C_{w+4-n}$ 
with 
$a_0=(\theta_0 - \theta_1 )/z^2$,\, 
$a_1= \{ \theta_1 (C_3+1)- \theta_0 - \theta_2 \}/z^2$,\,and $ a_2=(\theta_2 - \theta_1 )/z^2\, $.
With this choice $\Omega_{n,\,\beta}^{**} = \Omega_{n,\,\beta} - \Omega_{n,\,\beta}^{*}$ satisfies (\ref{Hconjformula}) and 
has $0=\theta_0^{**} = \theta_1^{**}=\theta_2^{**}$. In other words, $\Omega_{n,\,\beta}^{**}$ depends only on $h_j$ for $j \geq 3$.
A natural conjecture is that
$\Omega_{n,\,\beta}^{**}$ is zero or trivially defined as a weighted sum of the $h_j$.
What has been established is summarized in
%
%%%%%%%%%%%%%%%%%% Uniqueness of Equations satisfied by h_j %%%%%%%%%%%%%%%%%%%%%%%%
%
\begin{thm}
Any function, $\Omega_{n,\,\beta}$, that satisfies
$\sum_{j=0}^{n-1}   \,\theta_{j}\,h_j =  \Omega_{n,\,\beta}$, with $\theta_{j} \in \mathbb{Z}{[z, z^{-1}]}$
and $\beta \in B_n$,
has $\theta_j= \Omega_{n,\,\omega_j} $. Here $\omega_j$ is any $n$-braid of length $n-1$, writhe $n-1-2j$, 
and $\widehat{\omega_j}$ is the trivial knot. Further, we have
%$\Omega_{n,\,\beta}$ has the following properties:

\begin{enumerate}%[(i)]
\item when $\Omega_{n,\,\beta}$ is a link invariant,
$\Omega_{n,\,\beta}$ must be a multiple of the Conway polynomial for 
$\widehat{\beta}$; i.e., $\Omega_{n,\,\beta} = \theta_0~\nabla_{\widehat{\beta}}$, and
each $\theta_j$ equals $\theta_0$,
\item when $\Omega_{n,\,\beta}$ is not a link invariant, but depends only on $w-n$, we have
$\Omega_{n,\,\beta}= \theta_0 C_{3+w-n}/z - \theta_1C_{1+w-n}/z$\,,
and $\theta_j=  \{ \theta_0\,C_{-2j+2} - \theta_1\,C_{-2j} \}/z$\,,
\item when $\Omega_{n,\,\beta}$ is independent of the prior two classes,
we have $\theta_0, \theta_1, \theta_2$ are all zero. If $J$ is the minimal index for which $\theta_j \neq 0$,
we have $\theta_J h_J = \Omega_{J+1,\,\beta}$\,.
\end{enumerate}
\label{ThmUniqueLinEq}
\end{thm}
%$0=\theta_0=\theta_1=\theta_2$

%%%%%%%
A simple example shows why Cor.~\ref{HProperLcoeffrel} will be difficult to enhance.
Consider $\sigma_1^{1}$ and $\sigma_1^{-1}$ in $B_2$, 
for which $P_{\widehat{\sigma_1^{1}}}= P_{\widehat{\sigma_1^{-1}}}=1$ and $\mathtt{x}=\min \deg_v P=0$.
We have $C_{-3} h_{0,\,\sigma_1^{1} } = C_{\mathtt{x} +3}$,
and also $C_{-1} h_{1,\,\sigma_1^{-1} } = C_{\mathtt{x} +1}$.
This shows that when (\ref{HLaurentrel1}) is expressed using link invariants 
and the proper MFW polynomials, 
two different relations arise, though both braids represent the trivial knot.

This concludes the section of main results. It is time to apply these results
to see further properties of the \mbox{Homflypt} polynomial. Section~\ref{SubAppnBraids} will focus on 
general $n$-braid links, while
Section~\ref{SubAppThreeBraids} will only discuss three-braid links.

%%%%%%%%%%%%%%%%%%%%%%%%%%%%%%%%%%%%%%%  n-BRAID LINKS  %%%%%%%%%%%%%%%%%%%%%%%%%%%%%%%%%%%
%%%%%%%%%%%%%%%%%%%%%%%%%%%%%%%%%%%%%%%  n-BRAID LINKS  %%%%%%%%%%%%%%%%%%%%%%%%%%%%%%%%%%%
%%%%%%%%%%%%%%%%%%%%%%%%%%%%%%%%%%%%%%%  n-BRAID LINKS  %%%%%%%%%%%%%%%%%%%%%%%%%%%%%%%%%%%

\subsection{Applications to $n$-braid links}
\label{SubAppnBraids}

The following result applies to positive, and some non-alternating, braids.
\begin{cor}
Let $\beta$ be an $n$-braid whose first proper Laurent MFW polynomial (see Def.~\ref{DefMinMaxIndices}) attains 
the $z$-degree ($\zeta$) of $P_{\widehat{\beta}}$. 

The degree of $h_{ \mathtt{m} }$ satisfies $\zeta=\mathtt{x} = \min \deg_v P_{\widehat{\beta}}$ and $h_{ \mathtt{m} }$ is monic.

Also, $\mathtt{y} = \deg_v P \geq  \max\{(2(\mathtt{M}-\mathtt{m}) + 1-\mu(\widehat{\beta}),\,\mu(\widehat{\beta})-1 \}$.
\label{Cor1stProperHasMaxzdeg}
\end{cor}

\begin{proof}
By Cor.~\ref{HProperLcoeffrel} we may assume $b_{\min}(\beta) \geq 2$.
Choose $\kappa = 2 \mathtt{M}$ in $(\ref{HLaurentrel3})$, and observe that all the terms 
on the left side have odd negative subscripts to obtain 
$\sum_{j=\mathtt{m}}^{ \mathtt{M} }   \,C_{2 (\mathtt{M} -j) +1 }\,h_j =  C_{2 \mathtt{M} +2+w-n}$.
The unique highest degree term on the left is for $j=\mathtt{m}$, so its degree must be that of
the (non-zero) right side, so $\mathtt{y} +1=2 \mathtt{M} +2+w-n \neq 0$.
When $\mathtt{y} +1 >0$, the right side, $C_{2 \mathtt{M} +2+w-n}$, is monic, so $h_{ \mathtt{m} }$ is also monic; also
$\deg h_{ \mathtt{m} } = (2 \mathtt{M} +1+w-n) - 2 (\mathtt{M} - \mathtt{m})=\min \deg_v P_{\widehat{\beta}}
= \mathtt{x}$. %$2m+1+w-n$.
As $1-\mu(\widehat{\beta}) \leq \zeta =  \mathtt{x}$, we have 
$2(\mathtt{M}- \mathtt{m}) + 1-\mu(\widehat{\beta}) \leq \mathtt{y} $.
Apply
Lemma~\ref{KawauchiLemma}, with $s=0$, to see that 
the $v$-span ($\mathtt{y}-\mathtt{x}$) must be at least $2(\mu(\widehat{\beta}) -1)$, 
so $\mathtt{y} \geq \mathtt{x} + 2(\mu(\widehat{\beta}) -1) \geq 1-\mu(\widehat{\beta}) + 2(\mu(\widehat{\beta}) -1)$,
i.e. $ \mathtt{y} \geq \mu(\widehat{\beta}) - 1$.

When $\mathtt{y} +1 < 0$, we have $\deg h_{ \mathtt{m} } = -(2 \mathtt{M} +2+w-n)-1 - 2 (\mathtt{M} - \mathtt{m})$, i.e.
%$2m+ n-w-3-4 \mathtt{M}$, or 
$-\mathtt{y}-2 b_{\min}(\beta)$.
In order for $\deg h_{ \mathtt{m} } = \deg_z P_{\widehat{\beta}}$ when $b_{\min}(\beta) = 2$,
Cor.~\ref{HProperLcoeffrel} requires $\mathtt{x} \geq -1$, i.e. $\mathtt{y} \geq 1$,
so we assume $b_{\min}(\beta) \geq 3$ below.

The remainder of the proof that $\mathtt{y} > -1$ depends on Cor.~\ref{HProperLcoeffrel} (\ref{HProperLcoeffrel0});
note that $\nabla_{\widehat{\beta}}=0$ or $\deg \nabla_{\widehat{\beta}} \leq \deg h_{ \mathtt{m} }$.
When $\mathtt{y} \leq -2 $, we have $\deg C_{\mathtt{x}+3} = |\mathtt{x}| - 4$, so
$\deg z^{-2} C_{\mathtt{x}+3} = |\mathtt{x}| - 6 > |\mathtt{y}| - 2 b_{\min}(\beta) = \deg h_{ \mathtt{m} }$.
Thus there must be some $j \geq 3$ for which
$\deg  (C_{2j-3}-1) h_{\mathtt{m}+j} \geq \deg C_{\mathtt{x}+3} $, i.e.
$2j-4+ \deg h_{\mathtt{m}+j} \geq |\mathtt{x}| - 4 $.
As $|\mathtt{x}|=|\mathtt{y}| +  2( b_{\min}(\beta)-1)$, we need $2j + \deg h_{\mathtt{m}+j} \geq |\mathtt{y}| +  2( b_{\min}(\beta)-1)$
for some $j \in \left[3, \mathtt{M}- \mathtt{m} \right]$ and this is impossible 
if $\deg h_{\mathtt{m}+j} \leq \deg h_{ \mathtt{m} }$ %$= |\mathtt{y}| - 2 b_{\min}(\beta)$.
\end{proof}

The following is a restatement of Cor.~\ref{Cor1stProperHasMaxzdeg} for the mirror image braid.

\begin{cor}
Let $\beta$ be an $n$-braid whose last proper Laurent MFW polynomial attains 
the $z$-degree ($\zeta$)  of $P_{\widehat{\beta}}$\,. 

The degree of $h_{ \mathtt{M} }$ satisfies  $\zeta= - \deg_v P_{\widehat{\beta}}$ and $\epsilon_{w-n+1} h_{ \mathtt{M} }$ is monic. 

Also, $\mathtt{x} = \min \deg_v P_{\widehat{\beta}} \leq  \min\{\mu(\widehat{\beta}) -1 -2( \mathtt{M} -  \mathtt{m} ),\,1-\mu(\widehat{\beta}) \}$.
\label{CorLastProperHasMaxzdeg}
\end{cor}

We can now see that the only braids for which the 
first and last proper MFW polynomials (see Def.~\ref{DefMinMaxIndices}) both have the
$z$-degree of $P_{\widehat{\beta}}$ satisfy $P_{\widehat{\beta}} = P_{O_k}$ 
($O_k$ is the trivial link of $k$ components);
the link, $\widehat{\beta}$, might be non-trivial.

\begin{cor}
If $\beta$ is an $n$-braid whose first and last proper Laurent MFW polynomials both 
have the $z$-degree of $P_{\widehat{\beta}}$,
we have 
\begin{eqnarray*}
P_{\widehat{\beta}} &=& P_{O_{b_{\min}(\beta)}} = [(1-v^2)/vz]^{b_{\min}(\beta)-1}\,.
\end{eqnarray*}
\label{CorTrivialLinksUniqueness}
\end{cor}

\begin{proof}
By Cors.~\ref{Cor1stProperHasMaxzdeg} and \ref{CorLastProperHasMaxzdeg}
we have $\min \deg_v P_{\widehat{\beta}} = 
- \deg_v P_{\widehat{\beta}}$, so
$\min \deg_v P_{\widehat{\beta}}$ is $\mathtt{m} - \mathtt{M}$,
that is also $\deg_z P_{\widehat{\beta}}$ by Cor.~\ref{Cor1stProperHasMaxzdeg}.
Also $w - n + 1 + \mathtt{m} + \mathtt{M}=0$.

However, the minimum $z$-degree of $P_{\widehat{\beta}}$, i.e.  $1 - \mu(\widehat{\beta})$,
must not exceed the $z$-degree, i.e. $\mathtt{m} - \mathtt{M}$, so $\mu(\widehat{\beta}) -1 \geq  \mathtt{M}  - \mathtt{m} $.
By Lemma~\ref{KawauchiLemma}, with $s=0$, we also know $ \mathtt{M}  - \mathtt{m}  \geq \mu(\widehat{\beta})-1$,
so that $\mu(\widehat{\beta}) -1 =  \mathtt{M}  - \mathtt{m} $.
Since the minimum and maximum $z$-degrees are equal, each $h_j = a_j z^{\mathtt{m} - \mathtt{M}}$, for some integer $a_j$.

Choose $\kappa = 2 \mathtt{M}$ in $(\ref{HLaurentrel3})$ and substitute $w-n = - \mathtt{m} - \mathtt{M} -1$ to obtain \\
$\sum_{j= \mathtt{m}}^{\mathtt{M}}   \,C_{2 \mathtt{M} -2j +1 }\, a_j z^{m - \mathtt{M}} =  C_{b_{\min}(\beta)}$, i.e. 
$\sum_{j= \mathtt{m}}^{\mathtt{M}}   \,C_{2 \mathtt{M} -2j +1 }\, a_j =  z^{ \mathtt{M}  - \mathtt{m} } C_{b_{\min}(\beta)}$,
using $b_{\min}(\beta) =  \mathtt{M}  - \mathtt{m} +1$.
Now apply (\ref{Conwayzpowers}) to $z^{ \mathtt{M}  - \mathtt{m} }\,C_{b_{\min}(\beta)}$ to obtain

\begin{eqnarray*}
\sum_{j = \mathtt{m}}^{\mathtt{M}}   \,C_{2 (\mathtt{M} -j) +1 }\, a_j &=& 
\sum_{j=0}^{ \mathtt{M}  - \mathtt{m} }
\left( \begin{array} 
{cc}  \mathtt{M}  - \mathtt{m}  \\ j
\end{array} \right) 
\epsilon_{j}\,C_{2( \mathtt{M}  - \mathtt{m} )+1-2j}\,.
\end{eqnarray*}

This implies $a_{\mathtt{m}+j} = \epsilon_{j}
\left( \begin{array} 
{cc}  \mathtt{M}  - \mathtt{m}  \\ j
\end{array} \right)$,
for $j \in [0, \mathtt{M}  - \mathtt{m} ]$.
Eq.~\ref{HLaurentform} tells us
$P_{\widehat{\beta}} = v^{-\mathtt{m} - \mathtt{M}} 
\sum_{j=0}^{ \mathtt{M}  - \mathtt{m} } a_{\mathtt{m}+j} z^{\mathtt{m} - \mathtt{M}} v^{2(\mathtt{m}+j)}$, 
i.e. $(vz)^{\mathtt{m} - \mathtt{M}} (1-v^2)^{ \mathtt{M}  - \mathtt{m} }$.
\end{proof}

Equation~\ref{HLaurentrel3} may also be used to show that 
when there is a (unique) Laurent MFW polynomial, 
say $h_r$, that achieves the $z$-degree of $P_{\widehat{\beta}}$, 
and the $v$-span is positive,
there must be some proper Laurent MFW polynomial whose
degree is within twice its index distance from $r$.
 
\begin{cor}
If $\beta$ is an $n$-braid for which $\deg h_r = \deg_z P_{\widehat{\beta}}$,
and $b_{\min}(\beta) > 1$,
there is an index $s \neq r$ for which $h_s \neq 0$ and
$$\deg h_s \geq  \deg h_r - 2*|r-s|\,.$$
\label{CorMaxzdegreeIndexShadow}
\end{cor}

\begin{proof}
Choose $\kappa = 2 \mathtt{M}$ in $(\ref{HLaurentrel3})$ to obtain
$\sum_{j= \mathtt{m} }^{ \mathtt{M} }   \,C_{2 \mathtt{M} -2j +1 }\, h_j  =  C_{\mathtt{y} +1}$
(recall $\mathtt{y} = \deg_v P_{\widehat{\beta}}$).
If there is no $s \in \left [ \mathtt{m}, r \right)$ satisfying the assertion, 
$C_{2 \mathtt{M} -2r +1 }\, h_r$ has degree higher than any other $C_{2 \mathtt{M} -2j +1 }\, h_j$ and $\mathtt{y} \neq -1$.
Thus the degree of $h_r$ is the degree of $C_{\mathtt{y} +1}$ less $2(\mathtt{M} - r)$, i.e.
$|\mathtt{y} +1|-1-2( \mathtt{M} -r)$.

Similarly, the choice $\kappa = 2 \mathtt{m} -2$ in $(\ref{HLaurentrel3})$ gives us
$\sum_{j= \mathtt{m} }^{ \mathtt{M} }   \,C_{2 j -2 \mathtt{m} +1 }\, h_j  =  C_{\mathtt{x}-1}$
(recall $\mathtt{x} = \min \deg_v P_{\widehat{\beta}}$).
If there is no $s \in \left ( r, \mathtt{M} \right ]$ satisfying the assertion,
we see $\mathtt{x} \neq 1$ and the degree of $h_r$ is $|\mathtt{x} -1|-1- 2(r- \mathtt{m} )$.

If the assertion is false, a comparison of the two expressions for $\deg h_r$ shows
that $\mathtt{x} -1 < 0 <\mathtt{y} +1$ and $\mathtt{x}+\mathtt{y}+4r = 2(\mathtt{M} + \mathtt{m})$, so
%$w = n-2r-1$
$w - n+1+2r=0$, and the degree of $h_r$ is zero. This implies $\mu(\widehat{\beta} )$ is odd and $\geq 3$
(if $\mu(\widehat{\beta}) = 1$, we would have $P_{\widehat{\beta}}=1$ by Cor.~\ref{CorTrivialLinksUniqueness}).
%Since $h_r$ has degree higher than that of any other $h_j$, 
Apply Lemma~\ref{KawauchiLemma}, with $2s=\mu(\widehat{\beta} )-1$, to see 
that the $v$-span of $P_0(v)$ is at least $\mu(\widehat{\beta} )-1 \geq 2$, 
a contradiction.
\end{proof}

\begin{cor}
Suppose $\beta$ is an $n$-braid with 
$\zeta = \deg_z P_{\widehat{\beta}} > 1-\mu(\widehat{\beta})$.

When $\zeta > \max\{ \deg h_{ \mathtt{m} },\, \deg h_{ \mathtt{M} }\}$, 
we have $\zeta > \max\{\mathtt{x}, - \mathtt{y} \}$.

\label{CorvspanEndsNotMaxzDeg}
\end{cor}
%= C_{x+2( \mathtt{M}  - \mathtt{m} )+1}
\begin{proof}
Use $\kappa = 2 \mathtt{M}$ in
Eq.~\ref{HLaurentrel3} to see $\sum_{j= \mathtt{m} }^{ \mathtt{M} } \,C_{2(\mathtt{M} -j)+1}\, h_{j}  = C_{\mathtt{y} +1}$.
This implies $\deg C_{\mathtt{y}+1} \leq 2( \mathtt{M}  - \mathtt{m} ) + \zeta -2$.
For $\mathtt{y} \geq 0$ we have $\deg C_{\mathtt{y}+1} = \mathtt{y} = \mathtt{x} + 2( \mathtt{M}  - \mathtt{m} )$, so
$\mathtt{x} \leq \zeta -2$. 
Thus the result is true for $\{(\mathtt{x},\mathtt{y}): \mathtt{x} \geq - \mathtt{y} \mbox{ and } \mathtt{y} \geq 0 \}$.

For $\kappa=2 \mathtt{m} -2$ in Eq.~\ref{HLaurentrel3} we have
$\sum_{j= \mathtt{m} }^{ \mathtt{M} } \,C_{2(j- \mathtt{m})+1}\, h_{j}  = C_{\mathtt{x} -1}$. 
For $\mathtt{x} \leq 0$, we have 
$\deg C_{\mathtt{x} -1} = -\mathtt{x} \leq 2( \mathtt{M}  - \mathtt{m} ) + \zeta -2$, i.e.  $-\mathtt{y} < \zeta$.
We need only show $\zeta > \mathtt{x}$ for valid pairs $(\mathtt{x},\mathtt{y})$ with $\mathtt{x} \leq 0$ and $\mathtt{y} < -\mathtt{x}$.
Consider the value $\zeta - \mathtt{x}$ for a pair $(\mathtt{x},\mathtt{y})$ 
with $\mathtt{x} = \mathtt{y} - 2( \mathtt{M}  - \mathtt{m} )$ and
$\mathtt{y} <  \mathtt{M}  - \mathtt{m} $\,;
these are all the remaining cases.
By Lemma~\ref{KawauchiLemma} with $s=0$ we know $\mu(\widehat{\beta}) \leq b_{\min}(\beta)$, so we have
$\zeta - \mathtt{x} > 1-\mu(\widehat{\beta}) - \mathtt{x} \geq (\mathtt{m}  - \mathtt{M} ) - [\mathtt{y}-2( \mathtt{M}  - \mathtt{m} )]$\,, 
i.e. $\zeta - \mathtt{x} >  \mathtt{M}  - \mathtt{m}  - \mathtt{y} > 0$\,.
\end{proof}

%%%%%%%%%%%%%%  b_{\min} =3 THM VERSION %%%%%%%%%%%%%%%%%%%%%%%%%%%%%%%%%%%%%%
%%%%%%%%%%%%%%  b_{\min} =3 THM VERSION %%%%%%%%%%%%%%%%%%%%%%%%%%%%%%%%%%%%%%
%%%%%%%%%%%%%%  b_{\min} =3 THM VERSION %%%%%%%%%%%%%%%%%%%%%%%%%%%%%%%%%%%%%%

\begin{thm}
Suppose $\beta$ is an $n$-braid with 
$\zeta = \deg_z P_{\widehat{\beta}} > 1-\mu(\widehat{\beta})$ and $b_{\min}(\beta)=3$.
For $j=\mathtt{m}, \mathtt{m}+1, \mathtt{m}+2$, write $h_j$   
as in \ref{Thmvspan4hm}, \ref{Thmvspan4hmplus1}, \ref{Thmvspan4hmplus2} below and
define $a_k, b_k, c_k =0$ for $k \not \in [1-\mu(\widehat{\beta}), \zeta]$.
It follows that $\mu(\widehat{\beta}) \leq 3$ and $\zeta \geq 0$.

%Also $\zeta= \mathtt{x}$, or $\zeta= -\mathtt{y} > \mathtt{x}$, or $\zeta > \max\{\mathtt{x} ,-\mathtt{y} \}$.
$P_{\widehat{\beta}}$ is determined by $\zeta$, 
$\mathtt{x}$, any proper $h_J$, plus $J - \mathtt{m}$,
as described below. 
For $\mu(\widehat{\beta})=3$, we have $P_{-2}(v) = v^{ \mathtt{x} } (1-v^2)^2$ 
(see Def.~\ref{DefPk} for $P_{k}(v)$).

\begin{enumerate}%[(i)]
\item $h_{ \mathtt{m} }= a_{\zeta} z^{\zeta} + a_{\zeta-2} z^{\zeta-2} + \cdots + a_{1- \mu(\widehat{\beta})  } z^{1- \mu(\widehat{\beta}) }$\,,
\label{Thmvspan4hm}
\item $h_{m+1}= b_{\zeta} z^{\zeta} + b_{\zeta-2} z^{\zeta-2} + \cdots + b_{1- \mu(\widehat{\beta})  } z^{1- \mu(\widehat{\beta}) }$\,,
\label{Thmvspan4hmplus1}
\item $h_{m+2}= c_{\zeta} z^{\zeta} + c_{\zeta-2} z^{\zeta-2} + \cdots + c_{1- \mu(\widehat{\beta})  } z^{1- \mu(\widehat{\beta}) }$\,.
\label{Thmvspan4hmplus2}
\end{enumerate}

%%%%%%%%%%%%%%%%%%%%%%%%%%%% Case I %%%%%%%%%%%%%%%%%%%%%%%%%%%%%%%%%%%%%%%

\noindent
\textbf{Case I} ($\zeta \geq \mathtt{x} \geq -2$\,):
We have $c_{\zeta}=0$. If $c_{\zeta-2}=0$ we have
$\zeta=\mathtt{x}$, $a_{\zeta}=1$,  $a_{\zeta-2}=\zeta$, and $b_{\zeta}=0$.
$P_{-1}(v) = v^{ \mathtt{x} } \{ 1+c_{-1}  - (1+ 2c_{-1})v^2 + c_{-1} v^4 \}$ if $\mu(\widehat{\beta})=2$.
\begin{eqnarray} 
%a_{\zeta-2k} &=& c_{\zeta-2k} \,,\mbox{for } j \in \left[0, \frac{\zeta- \mathtt{x} -2}{2} \right]\,, \\
%
a_{\zeta-2j} &=&
\left( \begin{array} 
{cc}  \frac{\zeta + \mathtt{x} }{2} +1 -j \\ j - \frac{\zeta- \mathtt{x} }{2}
\end{array} \right)   
+ c_{\zeta-2j}\,,
%\mbox{for } j \in \left[ \frac{\zeta- \mathtt{x} }{2} ,\,
%1+\lfloor \frac {\zeta}{2} \rfloor \right]\,,  
\label{Thmvspan4CaseABCaj}  \\
%%%%%%%%% Thmvspan4CaseABCaj
%b_{\zeta-2j} &=& - 2 c_{\zeta-2j} - c_{\zeta-2-2j}\,,\mbox{for } 
%j \in \left[0, \frac{\zeta- \mathtt{x} -2}{2} \right]\,, \nonumber  \\
%%%%%%%%%%%%
b_{\zeta-2j} &=&
- \left( \begin{array} 
{cc}  \frac{\zeta + \mathtt{x} }{2}  -j \\ j-1 - \frac{\zeta- \mathtt{x} }{2}
\end{array} \right)
-  2 c_{\zeta-2j} - c_{\zeta-2-2j}\,, \mbox{ for }  \mathtt{x} \geq 0 \,, 
\label{Thmvspan4CaseABCbj} \\
b_{\zeta-2j} &=& \epsilon_{\mathtt{x}}  \delta_{j,\, \lfloor \frac {1+\zeta}{2} \rfloor } -  2 c_{\zeta-2j} - c_{\zeta-2-2j} \,, 
\mbox{ for } \mathtt{x} \in \{ -1, -2 \} \,.
\label{Thmvspan4CaseABCbjxisneg1}
\end{eqnarray}

%%%%%%%%%%%%%%%%%%%%%%%%%%%% Case II %%%%%%%%%%%%%%%%%%%%%%%%%%%%%%%%%%%%%%%

\noindent
\textbf{Case II} ($\zeta \geq \mathtt{-y} \geq -2$): this class is the mirror of braids in Case I.
\begin{eqnarray} 
b_{\zeta-2j} &=&
-\epsilon_{\mathtt{y}} \left( \begin{array} 
{cc}  \frac{\zeta - \mathtt{y} }{2}  -j \\ j-1 - \frac{\zeta+ \mathtt{y} }{2}
\end{array} \right)
-  2 a_{\zeta-2j} - a_{\zeta-2-2j}\,, \mbox{ for }  \mathtt{y} \leq 0 \,, 
\label{Thmvspan4Mirrorbj} \\
% following formula does not need  \epsilon_{\mathtt{y}} in front of \delta !!!
b_{\zeta-2j} &=&  \delta_{j,\, \lfloor \frac {1+\zeta}{2} \rfloor } -  2 a_{\zeta-2j} - a_{\zeta-2-2j} \,, 
\mbox{ for } \mathtt{y} \in \{ 1, 2 \}\,, 
\label{Thmvspan4Mirrorbjxisneg1}  \\
%%%%%%%%%%
%%%%%%%%%%
c_{\zeta-2j} &=& \epsilon_{\mathtt{y}}
\left( \begin{array} 
{cc}  \frac{\zeta - \mathtt{y} }{2} +1 -j \\ j - \frac{\zeta + \mathtt{y} }{2}
\end{array} \right)   
+ a_{\zeta-2j}\,.
\label{Thmvspan4Mirroraj}
\end{eqnarray}

Cases I and II cover all possibilities. The relations in each of the cases
may be manipulated so that the independent coefficients are $a_k$ or $b_k$ or $c_k$.
\label{Thmvspan4}
\end{thm}
%%%%%%%%%%%%%%%%%%%%%%%%%%%% PROOF   %%%%%%%%%%%%%%%%%%%%%%%%%%%%%%%%%%%%%%%
\begin{proof}
As $b_{\min}(\beta)=3$, Lemma~\ref{KawauchiLemma} ($s=0$) gives $\mu(\widehat{\beta}) \leq 3$. 
If $s \in [0, \mu(\widehat{\beta})-2]$, (p.133, ~\cite{51}) tells us
$P_{2s+1-\mu(\widehat{\beta})}(1)=0$.
If $\mu(\widehat{\beta})=3$, we see $P_{-2}(v)=c_{-2}v^{\mathtt{x}} (1-v^2)^2$; now
evaluate (\ref{Thmvspan4CaseABCaj}, \ref{Thmvspan4CaseABCbj}, \ref{Thmvspan4CaseABCbjxisneg1}) at 
$j=\frac{\zeta}{2}$ to see $c_{-2}=1$ as $P_{0}(1)=0$.

Cor.~\ref{Cor1stProperHasMaxzdeg} addresses the only case where $\zeta= \mathtt{x}$ (Case I);
Cors.~\ref{CorLastProperHasMaxzdeg} and \ref{CorvspanEndsNotMaxzDeg} have $\zeta> \mathtt{x}$ 
(for Cor.~\ref{CorLastProperHasMaxzdeg} this is due to the assumption
$\zeta > 1-\mu(\widehat{\beta})$). 
Cor.~\ref{CorLastProperHasMaxzdeg} addresses the only case where $\zeta= - \mathtt{y}$ (Case II) as
Cor.~\ref{CorvspanEndsNotMaxzDeg} has $\zeta> \mathtt{x}, -\mathtt{y}$.
Cor.~\ref{CorTrivialLinksUniqueness} is excluded by 
the assumption that $\zeta > 1-\mu(\widehat{\beta})$.
These corollaries cover all cases, i.e. $\zeta = \deg h_{ \mathtt{m} }$ or $\zeta = \deg h_{ \mathtt{M} }$
or $\zeta > \max\{ \deg h_{ \mathtt{m} }, \deg h_{ \mathtt{M} }   \}$.
Since $b_{\min}(\beta)=3$, we have $\mathtt{y} = \mathtt{x} +4$, 
so $\mathtt{x} \geq -2$ is equivalent to $-\mathtt{y} \leq -2$.
Thus Cases I and II cover all possibilities and overlap at $\mathtt{x} = -2 =-\mathtt{y}$.
In Case I, we have $c_{\zeta}=0$. If $c_{\zeta-2}=0$, then $\zeta \neq 0$ and
(\ref{Thmvspan4CaseABCaj}, \ref{Thmvspan4CaseABCbj}, \ref{Thmvspan4CaseABCbjxisneg1}) imply
$\zeta=\mathtt{x}$ to avoid $a_{\zeta}=b_{\zeta}=0$; now 
Eqs.~\ref{Thmvspan4CaseABCaj}, \ref{Thmvspan4CaseABCbj} imply $a_{\zeta}=1$, $a_{\zeta-2}=\zeta$ and $b_{\zeta}=0$.

%%%%%%%%%%%%%%%%%%%%%%%%%%%%% Case I PROOF  %%%%%%%%%%%%%%%%%%%%%%%%%%%%%%%%%%%%%%

In Case I use (\ref{HLaurentrel3}) with $\kappa = 2 \mathtt{m} +2$  to see 
$\sum_{j=0}^{2}   \,C_{-3+2 j}\, h_{\mathtt{m} +j}  =  C_{\mathtt{x} +3}$ 
and use $\kappa = 2 \mathtt{m}$ to see $\sum_{j=0}^{2}   \,C_{-1+2 j}\, h_{\mathtt{m} +j}  =  C_{\mathtt{x} +1}$.
The coefficients of $z^{\zeta-2j}$ in these expressions are
$a_{\zeta-2j}+ a_{\zeta-2-2j} + b_{\zeta-2j}+ c_{\zeta-2j}$ and
$a_{\zeta-2j} + b_{\zeta-2j}+ c_{\zeta-2j} + c_{\zeta-2-2j}$ respectively.
With $\xi = (\zeta- \mathtt{x} -2)/2$, expand $C_{\mathtt{x} +3}$ using
(\ref{Ctorusstdform}) to see
%
% C_{x+3} below
\begin{eqnarray}
0 =  a_{\zeta-2j}+ a_{\zeta-2-2j} + b_{\zeta-2j}+ c_{\zeta-2j}\,,\mbox{ for } 
%j < \xi  \,,\mbox{ or } j > \lfloor \frac {\zeta}{2} \rfloor \,,\mbox{ and } &&
j < \xi  \,,\mbox{ or } j > \lfloor \zeta/2 \rfloor \,,\mbox{ and } &&
\label{CaseABCCxplus3zero} 
\end{eqnarray}
\begin{eqnarray}
\left( \begin{array} 
{cc} \frac{\zeta+ \mathtt{x} +2}{2}-j \\ j -  \frac{\zeta- \mathtt{x} -2}{2}
\end{array} \right)   
= a_{\zeta-2j}+ a_{\zeta-2-2j} + b_{\zeta-2j}+ c_{\zeta-2j}\,,
%\mbox{ for } j \in \left[ \xi ,\,\lfloor \frac {\zeta}{2} \rfloor \right].&&
\mbox{ for } j \in \left[ \xi ,\,\lfloor \zeta/2 \rfloor \right].&&
\label{CaseABCCxplus3nonzero}
\end{eqnarray}
%
%% C_{x+1} below
%

We now assume $\mathtt{x} \geq 0$ until otherwise specified.
The expression for $C_{\mathtt{x} +1}$ given by (\ref{Ctorusstdform}) yields the following:
\begin{eqnarray}
0 =  a_{\zeta-2j} + b_{\zeta-2j}+ c_{\zeta-2j} + c_{\zeta-2-2j}\,,\mbox{ for } 
j \leq \xi \,,
%\mbox{ or } j > \lfloor \frac {\zeta}{2} \rfloor \,, \mbox{ and }&&
\mbox{ or } j > \lfloor \zeta/2 \rfloor \,, \mbox{ and }&& 
\label{CaseABCCxplus1zero} 
\end{eqnarray}
\begin{eqnarray} 
\left( \begin{array} 
{cc}  \frac{\zeta + \mathtt{x} }{2} -j \\ j - \frac{\zeta- \mathtt{x} }{2}
\end{array} \right)   
= a_{\zeta-2j} + b_{\zeta-2j}+ c_{\zeta-2j}+ c_{\zeta-2-2j}\,,
%\mbox{ for } j \in \left[ 1+\xi ,\,\lfloor \frac {\zeta}{2} \rfloor \right]. &&
\mbox{ for } j \in \left[ 1+\xi ,\,\lfloor \zeta/2 \rfloor \right]. &&
\label{CaseABCCxplus1nonzero}
\end{eqnarray}

Compare (\ref{CaseABCCxplus3zero}, \ref{CaseABCCxplus1zero}) to see $a_{\zeta-2-2j} = c_{\zeta-2-2j}$ for
$j < \xi $ or $j > \lfloor \frac {\zeta}{2} \rfloor$.
%$\mathbf{ a_{\zeta-2k} = c_{\zeta-2k} }$ for $k \in \left[1, \frac{\zeta- \mathtt{x} -2}{2} \right]$; actually works for $k=0$ also.
For $j=\xi$, (\ref{CaseABCCxplus3nonzero}, \ref{CaseABCCxplus1zero}) imply
$ 1 =  a_{\mathtt{x}} - c_{\mathtt{x}} $.

Subtract (\ref{CaseABCCxplus1nonzero}) from (\ref{CaseABCCxplus3nonzero}) and use (\ref{Pascaltriangle}) to obtain, 
for $j \in \left[ 1+\xi ,\,\lfloor \frac {\zeta}{2} \rfloor \right]$:
\begin{eqnarray} 
\left( \begin{array} 
{cc}  \frac{\zeta + \mathtt{x} }{2} -j \\ j +1 - \frac{\zeta- \mathtt{x} }{2}
\end{array} \right)   
&=&  a_{\zeta-2-2j} - c_{\zeta-2-2j}\,.
\label{temp1}
\end{eqnarray}
Eq.~\ref{temp1} is valid for any $j \leq \lfloor \frac {\zeta}{2} \rfloor $ by the prior paragraph. 
For $j=\lfloor \frac {\zeta}{2} \rfloor$ and odd $\zeta$
we see $1=a_{-1} - c_{-1}$.
For $j > \lfloor \frac {\zeta}{2} \rfloor$ all terms in (\ref{temp1}) are zero.
Eq.~\ref{temp1} is (\ref{Thmvspan4CaseABCaj}).
%Replace $j+1$ by $j$ in $\ref{temp1}$ to obtain (\ref{Thmvspan4CaseABCaj}).

For $j \leq \xi $, Eq.~\ref{Thmvspan4CaseABCbj} follows from (\ref{Thmvspan4CaseABCaj}, \ref{CaseABCCxplus1zero}).
For $j =1+ \xi $, Eq.~\ref{Thmvspan4CaseABCbj} follows from (\ref{CaseABCCxplus1nonzero})
%$1= a_{\mathtt{x}} + b_{\mathtt{x}}+ c_{\mathtt{x}}+ c_{\mathtt{x} -2}$,
and $1 =  a_{\mathtt{x}} - c_{\mathtt{x}} $ (see above or use Eq.~\ref{Thmvspan4CaseABCaj}).

Eqs.~\ref{Thmvspan4CaseABCaj}, \ref{CaseABCCxplus1nonzero} for
$j \in \left[ 1+\frac{\zeta- \mathtt{x} }{2}\,, \lfloor \frac {\zeta}{2} \rfloor \right]$ plus (\ref{Pascaltriangle}) imply
(\ref{Thmvspan4CaseABCbj}).
For $j > 1+\lfloor \frac {\zeta}{2} \rfloor$ all terms in (\ref{Thmvspan4CaseABCbj}) are zero.
For $j = 1+\lfloor \frac {\zeta}{2} \rfloor$ and odd $\zeta$, (\ref{Thmvspan4CaseABCbj}) yields
$b_{-1} = -1 - 2c_{-1}$, as required by (\ref{Thmvspan4CaseABCaj}, \ref{CaseABCCxplus3zero}).
For $j = 1+\lfloor \frac {\zeta}{2} \rfloor$ and even $\zeta$, (\ref{Thmvspan4CaseABCbj}) yields
$b_{-2} =   - 2c_{-2}$, as required by (\ref{Thmvspan4CaseABCaj}, \ref{CaseABCCxplus3zero}). 
This establishes (\ref{Thmvspan4CaseABCbj}) for $\mathtt{x} \geq 0$.

In case $\mathtt{x} = -1$ we see
$a_{\zeta-2j}+ a_{\zeta-2-2j} + b_{\zeta-2j}+ c_{\zeta-2j} = \delta_{ j,\,\lfloor \frac {\zeta}{2} \rfloor }$ 
and $a_{\zeta-2j} + b_{\zeta-2j}+ c_{\zeta-2j}+ c_{\zeta-2-2j} = 0$ for all $j$.
Use $j = \lfloor \frac {\zeta}{2} \rfloor$ to see $a_{-1} - c_{-1} = 1$; otherwise
$a_{\zeta-2-2j} = c_{\zeta-2-2j}$. This establishes (\ref{Thmvspan4CaseABCaj});
% a_k=c_k for k \neq -1$
Eq.~\ref{Thmvspan4CaseABCbjxisneg1} follows immediately.
%As $b_{\zeta-2j} =  - a_{\zeta-2j} - c_{\zeta-2j} - c_{\zeta-2-2j}$, we are done.

In case $\mathtt{x} =-2 $ we have
$a_{\zeta-2j}+ a_{\zeta-2-2j} + b_{\zeta-2j}+ c_{\zeta-2j} = \delta_{ j,\, \frac {\zeta}{2} }$ and
$a_{\zeta-2j} + b_{\zeta-2j}+ c_{\zeta-2j}+ c_{\zeta-2-2j} = \delta_{ j,\, \frac {\zeta}{2} }$.
Hence $a_{\zeta-2j} = c_{\zeta-2j}$ for all $j$ (this is Eq.~\ref{Thmvspan4CaseABCaj}) and
$ b_{\zeta-2j} =  \delta_{ j,\, \frac {\zeta}{2} } - 2 c_{\zeta-2j} -c_{\zeta-2-2j}$
which establishes (\ref{Thmvspan4CaseABCbjxisneg1}).

In Case I, it is clear that $a_k$ can be made the independent coefficient. To see that this is also true
for $b_k$, observe that $b_{1- \mu(\widehat{\beta})}$ is a function of $c_{1- \mu(\widehat{\beta})}$ that can be inverted; 
in general $b_{\zeta-2j}$ is a function of $c_{\zeta-2j}$ and $c_{\zeta-2-2j}$, hence $c_{\zeta-2j}$ can be expressed
in terms of $\{b_k\}$ (similiarly for $a_{\zeta-2j}$). A similar argument applies to Case II.
\end{proof}

%%%%%%%%%%%%%%%%%%%%%%%%%%%%%%%%%%%%%%%  3-BRAID LINKS  %%%%%%%%%%%%%%%%%%%%%%%%%%%%%%%%%%%
%%%%%%%%%%%%%%%%%%%%%%%%%%%%%%%%%%%%%%%  3-BRAID LINKS  %%%%%%%%%%%%%%%%%%%%%%%%%%%%%%%%%%%
%%%%%%%%%%%%%%%%%%%%%%%%%%%%%%%%%%%%%%%  3-BRAID LINKS  %%%%%%%%%%%%%%%%%%%%%%%%%%%%%%%%%%%

\subsection{Applications to three-braid links}
\label{SubAppThreeBraids}

All link and braid references in this section are to three-braid links and $B_3$.
It is convenient to restate Thm.~\ref{Hcoeffrel} as Prop.~\ref{H3braidstdformprop}. 
For three-braid knots Eq.~8.4, p.~356, \cite{41}, is equivalent to
Prop.~\ref{H3braidstdformprop}, that is valid for all three-braid links.
H.~Murakami, corollary \cite{69}, shows that
$P_L(a,t^{\frac{1}{2}} - t^{-\frac{1}{2}})$ may be expressed in terms of weighted sums of $1$,
$\Delta_L(t)$, where the weights are certain
Laurent polynomials in $a, \sqrt t $ for diagrams with 3 Siefert circles.
This appears to be the logical basis of 
Prop.~\ref{H3braidstdformprop}.

\begin{prop} Three-braid links, $\widehat \beta$, with $\beta \in B_3$, satisfy (\ref{Hstdform}) with
\begin{enumerate}%[(i)]
\item $p_0 = C_{w+1} - \nabla_{\widehat{\beta}}$\,,
\item $p_{1} =  z^2\, \nabla_{\widehat{\beta}} - p_0 - p_2 $\,,
\item $p_2 = C_{w-1} - \nabla_{\widehat{\beta}}$\,.
\end{enumerate}

\label{H3braidstdformprop}
\end{prop}

%%%%%%%%%%%%%%%%%%%%%%%% 3-BRAID Homflypt %%%%%%%%%%%%%%%%%%%%%%%%%%%%%%%%%%%%%%%%%
An equivalent form is Lemma~\ref{Lemma3brHomfly}, with 
Props.~\ref{Writhe3brPlusAny} and \ref{HConway3br} as easy corollaries.
J.~Birman,  \cite{5}, first established Prop.~\ref{Writhe3brPlusAny}
by use of Burau representations of three-braids.
%showed that when 
%a pair of three-braids have the same writhe and their closures have the same
%Jones (or Alexander) polynomial, then their \mbox{Homflypt} polynomials are also equal,

\begin{lem}
The \mbox{Homflypt} polynomial for a link, $\widehat{\beta}$\,, with 
$\beta \in B_3$ is 
\begin{equation}
P_{\widehat{\beta}} = P_{T_w}\,P_{O_2} - \nabla_{\widehat{\beta}}~v^w\,(P_{O_3}-1)\,.
\label{H3braidorigstdform}
\end{equation}

$O_2$ and $O_3$ are the trivial links of two and three components, respectively, so
$P_{O_2} = (1-v^2)/vz$ and $P_{O_3} = (1-v^2)^2/(vz)^2$.  
$P_{T_w}$ is reflected in(\ref{Htorusstdform}).
\label{Lemma3brHomfly}
\end{lem}

\begin{prop} If $\beta,\, \gamma \in B_3$\,, have the same writhe,
and $\widehat \beta$ and $\widehat \gamma$ have the same
Conway, Jones, or Alexander polynomials, we have 
$P_{\widehat \beta} = P_{\widehat \gamma}$\,.
\label{Writhe3brPlusAny}
\end{prop}

\begin{prop} When three-braid words, $\beta$, $\gamma$, satisfy 
$w(\beta) \geq w(\gamma)$\,, we have $P_{\widehat \beta} = P_{\widehat \gamma}$
exactly when $\nabla_{\widehat \beta} = \nabla_{\widehat \gamma}$ and 
one of the following is true:
\begin{enumerate}%[(i)]
\item $w(\beta) = w(\gamma)$\,, or
\label{HConway3brEqualw}
\item $w(\beta) = w(\gamma) +2$\,, and $\nabla_{\widehat \beta} = C_{w(\beta)-1}$\, or
\label{HConway3brw2apart}
\item $w(\beta) =2$\,, $w(\gamma)=-2$\,, and $\nabla_{\widehat \beta} = 1$\,.
\label{HConway3brtrivial}
\end{enumerate}
\label{HConway3br}
\end{prop}

\begin{proof}
Assume that $P_{\widehat \beta} = P_{\widehat \gamma}$\,.
If (\ref{HConway3brEqualw}) is false, we have $w(\beta) > w(\gamma)$\,.
Eq.~\ref{Hstdform} implies $w(\beta) = w(\gamma)+ 2$ or 
$w(\beta) = w(\gamma)+ 4$; and  $p_{2,\, \beta} =0=p_{0,\, \gamma}$.
Prop.~\ref{H3braidstdformprop} implies 
$\nabla_{\widehat \beta} = C_{w(\beta)-1}$\, and $\nabla_{\widehat \gamma} = C_{w(\gamma)+1}$\,.
Since $C_{w(\beta)-1} = C_{w(\gamma)+1}$\,,
we have $w(\beta) = w(\gamma)+2$ or $w(\beta) = -w(\gamma)$\,, i.e.
(\ref{HConway3brw2apart}) or (\ref{HConway3brtrivial}).

Conversely, when $\nabla_{\widehat \beta} = \nabla_{\widehat \gamma}$\,,
the first condition implies $P_{\widehat \beta} = P_{\widehat \gamma}$\, by Prop.~\ref{Writhe3brPlusAny}.
Prop.~\ref{H3braidstdformprop} shows the other two cases also imply 
$P_{\widehat \beta} = P_{\widehat \gamma}$\,.
\end{proof}

%%%%%%%%%%%%%%%  Properties of the Homflypt polynomial for three-braid links %%%%%%%%%%%%%%%%%%%%%%%%%%%
Prop.~\ref{Conway3braugmented} and Prop.~\ref{H3braidstdformprop} imply the following 
for three-braid links:
\begin{prop}
When $\gamma \in B_3$ and $a >0$ we have:
\begin{eqnarray*}
P_{\widehat {\alpha_2^3 \gamma} } &=& v^6 P_{\widehat {\gamma} } + P_{T_{w(\gamma)+5}} - v^6 P_{T_{w(\gamma)+1}} \,, \\
P_{\widehat {\alpha_2^{3a} \gamma} } &=& 
v^{6a} P_{\widehat {\gamma} }+ \sum_{j=1}^a v^{6a-6j} P_{T_{w(\gamma)+6j-1}} - \sum_{j=1}^a v^{6+6a-6j} P_{T_{w(\gamma)+6j-5}}\,, \\
P_{\widehat {\alpha_2^{-3a} \gamma} } &=& 
v^{-6a} P_{\widehat {\gamma} }+ \sum_{j=1}^a v^{6j-6a} P_{T_{w(\gamma)-6j+1}} - \sum_{j=1}^a v^{6j-6a-6} P_{T_{w(\gamma)-6j+5}}\,.
\end{eqnarray*}
\label{Homflypt3braugmented}
\end{prop}
The Jones polynomial for three-braid links may be derived from (\ref{H3braidorigstdform}).
Eq.~\ref{V3brformula} matches the knot formula in Prop.~11.10, p.~366,  \cite{41}, 
but the presence of $\epsilon_w$ in (\ref{V3brformula}) accomodates all 
three-braid links. With $G = t^2+t+1$, we have
%There is a difference in the two-variable polynomials for $O_2$ 
%(and also $O_{2k}$, but not $O_{2k+1}$), so that 
%Ex.~6.6, p.~350 \cite{41} gives the value as $(t - t^{-1})/x$, while in
%K.~Murasugi's book, Ex.~11.3.2, p.~232 \cite{71}, and \cite{45}, it is $(1-v^2)/vz = (v^{-1}-v)/z$.
%This appears to be due to a typo in \cite{41}, as the substitutions in $X_L(q,\,\lambda)$ to obtain $P_L(t, x)$
%yield $(t^{-1} - t)/x$\,.
\begin{equation}
V_{\widehat{\beta}}(t) = t^{(w-2)/2} \{t^{w+1} + \epsilon_w G \} 
- G t^{w-1}\Delta_{\widehat{\beta}}(t)\,, \mbox{ for } \beta \in B_3\,.
\label{V3brformula}
\end{equation}

These results may be combined with the extensive analysis by K.~Murasugi, \cite{70}, to 
determine all cases when the $v$-span of the \mbox{Homflypt} polynomial for a three-braid link is $4$ (or $0, 2$).
Thm.~\ref{Thmvspan4} and Prop.~\ref{Conway3braugmented} may be used to identify possible values for 
$(\mathtt{x},\, \zeta, \{c_j\})$ in Thm.~\ref{Thmvspan4}, but
a complete solution requires a manageable expression for $\nabla_{\widehat{\eta}}$ when $\eta$ is
an alternating three-braid.

J.~Birman, \cite{5}, first identified examples when the $v$-span is $2$
and further showed the Jones polynomial is a stronger invariant than the 
Alexander polynomial.
Research in \cite{105} suggests the Jones polynomial may distinguish the same 
three-braid links as the \mbox{Homflypt} polynomial.  
This is confirmed by Prop.~\ref{VequivP3br} and 
answers Question~4.1, p.~18 \cite{105}. Prop.~\ref{VequivP3br} is a consequence of
Lemma~\ref{HJones3br} and the Murasugi classification (Prop.~\ref{Murasugi3brpartitioning}).

In \cite{70}, K.~Murasugi defines a 
collection of seven disjoint sets, $\Omega_i$, of three-braid words and shows (Prop.~2.1 p.~7)
that each three-braid word is conjugate to exactly one element in some $\Omega_i$.
As the numbering and definition of the Artin braid generators 
differs from that in the later text, \cite{71}, we present a similar partitioning,
$\Omega_i^*$, consistent with the latter.
The Alexander polynomial for a typical member may be derived by use of 
(\ref{DeltaAnyPowerCenter}) and is included below; (\ref{MaxDegForDeltaInOmega0to5})
follows by inspection.

\begin{prop}
Each three-braid is conjugate to a unique element in some $\Omega_i^*$.
The variables $e,E, r, e_k, E_k $ are positive integers while $d$ is any integer.
Set $G = t^2+t+1 = (t^3-1)/(t-1)$.
For $\beta \in \Omega_j^*$ with $j \not = 6$, 
\begin{eqnarray}
\deg \Delta_{\widehat{\beta}} &=& (w(\beta)-2)/2\,.
\label{MaxDegForDeltaInOmega0to5}
\end{eqnarray}

\begin{enumerate}%[(i)]
\item $\Omega_0^* = \{\alpha_2^{3d}  \}$,
with $\Delta_{\widehat{\alpha_2^{3d}} } =  t(t^{3d}-1)^2/(Gt^{3d})$,
\label{Murasugi3brOmega0}
\item $\Omega_1^* = \{\alpha_2^{3d+1}  \}$,
with $\Delta_{\widehat{\alpha_2^{3d+1}} } = t(t^{6d+2} + t^{3d+1}+ 1)/(Gt^{3d+1}) $,
\label{Murasugi3brOmega1}
\item $\Omega_2^* = \{\alpha_2^{3d+2}  \}$,
with $\Delta_{\widehat{\alpha_2^{3d+2}} } =  t(t^{6d+4} + t^{3d+2}+ 1)/(Gt^{3d+2})$,
\label{Murasugi3brOmega2}
\item $\Omega_3^* = \{\alpha_2^{3d+1} \sigma_2  \}$,
with $\Delta_{\widehat{\alpha_2^{3d+1} \sigma_2} } = t(t^{6d+3} - 1)/(Gt^{3d+1}\sqrt{t})$,
\label{Murasugi3brOmega3}
\item $\Omega_4^* = \{\alpha_2^{3d} \sigma_2^{-e}  \}$, 
with $\Delta_{\widehat{\alpha_2^{3d} \sigma_2^{-e}} } = t (t^{3d}-1) (t^{3d-e}-\epsilon_e) /(Gt^{3d} t^{-e/2}) $,
\label{Murasugi3brOmega4}
\item $\Omega_5^* = \{\alpha_2^{3d} \sigma_1^{E} \}$, 
with $\Delta_{\widehat{\alpha_2^{3d} \sigma_1^{E}} } = t (t^{3d}-1) (t^{3d+E}-\epsilon_E) /(Gt^{3d} t^{E/2})$,
\label{Murasugi3brOmega5}
\item $\Omega_6^* = \{\alpha_2^{3d} \eta:  \eta \in B_3, 
\mbox{with $\eta = \prod_{k=1}^{r} \sigma_2^{-e_k} \sigma_1^{E_k}$  }  \}$, \newline
with $\Delta_{\widehat{\alpha_2^{3d} \eta} } = 
\Delta_{\widehat \eta} +  t (t^{3d}-1) (t^{3d+w(\eta)}-\epsilon_{w(\eta)}) /(Gt^{3d} t^{w(\eta)/2})$.
\label{Murasugi3brOmega6}
\end{enumerate}
\label{Murasugi3brpartitioning}
\end{prop}

%%%%%%%%%%%%%%%%%%%%%%%%%%%%%%%%%%%% LEMMA %%%%%%%%%%%%%%%%%%%%%%%%%%%%%%%%%%%%%%
%%%%%%%%%%%%%%%%%%%%%%%%%%%%%%%%%%%% LEMMA %%%%%%%%%%%%%%%%%%%%%%%%%%%%%%%%%%%%%%
%%%%%%%%%%%%%%%%%%%%%%%%%%%%%%%%%%%% LEMMA %%%%%%%%%%%%%%%%%%%%%%%%%%%%%%%%%%%%%%

\begin{lem} 
Assume $\beta, \gamma$ are three-braids. Set $a=w(\beta)$ and $b=w(\gamma)$. 
Suppose $a >  \max(b, -1)$ and $V_{\widehat \beta} = V_{\widehat \gamma}$.
We have $\mu(\widehat\beta) \not= 3$ and $a = b + 2k$.

For $k=1$, we have $\nabla_{\widehat \beta} = C_{w(\beta) - 1}=C_{w(\gamma) + 1}=\nabla_{\widehat \gamma}$\,.

When $\widehat \beta$ is a knot with $k>1$, we have  
$k=2$ and $a \geq 2$.
If $a = 2$, we have $\nabla_{\widehat \beta} =1= \nabla_{\widehat \gamma}$\,.
For $a > 2$, we have $a = 8p+2$ and $b = 8p-2 \geq 6$ and
\begin{eqnarray}
%t^{\frac{a-2}{2} } \Delta_{\widehat{\beta}}(t) &=&  t^{8p} + \sum_{j=0}^{2p-1} \epsilon_{j} t^{4j}(1 - t+t^3),
%\label{HJones3brDeltafor3brknot}  \\
%%%
t^{\frac{b}{2} } \Delta_{\widehat{\gamma}}(t)
%&=& 1 -t^{2}  + t^{3}\sum_{j=0}^{2p-2} \epsilon_{j} t^{4j}(1 - t+t^3).\\
&=& (t^{2} -1)t^{8p-4}+ \sum_{j=0}^{2p-2} \epsilon_j t^{4j} (1-t^2+t^3)\,.
\label{HJones3brDeltafor3brknotlow} 
\end{eqnarray}

For $\mu(\widehat{ \beta}) = 2$ and $k>1$, we have $k$ odd, 
$a=k(4p+3)$ and $b \geq 3$ and 

\begin{eqnarray}
t^{\frac{b}{2} } \Delta_{\widehat{\gamma}}(t)
=t^{4kp} (t^{k} - 1) + \sum_{j=0}^{p-1} t^{4kj} \{-1+ t^{k} - \sum_{x=0}^{k-1} t^{3x+k} (t-t^2)  \}\,.&& 
% original version
%= - 1 + t^{k}
%+  t^{k+1}  \sum_{j=0}^{p-1} t^{4 k j} [t^{3k-1}( t^{k} -1 ) + \sum_{x=0}^{k-1} t^{3x}  (t-1) ]\,. &&
\label{HJones3brDeltafor3br2linklow}
\end{eqnarray}

\label{HJones3br}
\end{lem} 

\begin{proof}
Observe that $V_{\widehat \beta} = V_{\widehat \gamma}$ 
implies $\mu(\widehat \beta)=\mu(\widehat \gamma)$ since $V_L(1)=(-2)^{\mu(L) -1}$,
p.~368 \cite{41}.
%or Thm.8.4.2 \cite{54}.
As $w-n \equiv \mu \mbox{ mod }2$, we have $a=b+2k$, with $k > 0$.

Comparison of the Jones polynomials, (\ref{V3brformula}), for $\widehat \beta$ and $\widehat \gamma$ yields an expression, 
(\ref{V3brequalDelta1}), for $\Delta_{\widehat{\gamma}}(t)$ in terms of $\Delta_{\widehat{\beta}}(t)$;
recall that $G = t^2+t+1$.
\begin{eqnarray}
t^{b}\Delta_{\widehat{\gamma}}(t)  &=& 
(t^{3b/2} -t^{3a/2} )t/G
+ \epsilon_a (t^{b/2} -t^{a/2})
+t^{a}\Delta_{\widehat{\beta}}(t)\,. \label{V3brequalDelta1}
\end{eqnarray}

Now replace $t$ by $1/t$ and multiply both sides in (\ref{V3brequalDelta1}) by $\epsilon_a$ to obtain
\begin{eqnarray*}
\epsilon_a t^{-b} \Delta_{\widehat{\gamma}}(1/t)  &=& 
\epsilon_a (t^{-3b/2} -t^{-3a/2} )t/G
+  (t^{-b/2} -t^{-a/2})
+ \epsilon_at^{-a}\Delta_{\widehat{\beta}}(1/t)\,. 
\end{eqnarray*}

Use
$\Delta_{L}(1/t) = (-1)^{\mu(L)-1}\Delta_{L}(t)$ , p.~6 \cite{105}, 
and multiply both sides in the prior line 
first by $t^{2a}$ and then multiply the prior line
by $t^{2b}$ to obtain

\begin{eqnarray*}
t^{2a-b} \Delta_{\widehat{\gamma}}(t)  &=& 
\epsilon_a (t^{(4a-3b)/2} -t^{a/2} )t/G
+  (t^{(4a-b)/2} -t^{3a/2})
+  t^{a} \Delta_{\widehat{\beta}}(t)\,, \\
t^{b} \Delta_{\widehat{\gamma}}(t)  &=& 
\epsilon_a (t^{b/2} -t^{(4b-3a)/2} )t/G
+  (t^{3b/2} -t^{(-a+4b)/2})
+  t^{2b-a}\Delta_{\widehat{\beta}}(t)\,. 
\end{eqnarray*}

Use these two relations with (\ref{V3brequalDelta1}) to replace 
$t^{a} \Delta_{\widehat{\beta}}(t)$ and $t^{b} \Delta_{\widehat{\gamma}}(t)$.
This yields equations for $(t^{2a-b} - t^{b}) \Delta_{\widehat{\gamma}}(t)$
and $(t^{a}-  t^{2b-a})\Delta_{\widehat{\beta}}(t)$ in terms of $t$.
Use the relation $a=b+2k$ to obtain

\begin{eqnarray}
t^{b/2}(t^{4k} - 1) \Delta_{\widehat{\gamma}}(t)  =
(t^{b}+\epsilon_b t^{k} ) (t^{3k}- 1 )t/G 
+ (t^{b+3k} +\epsilon_b)(t^{k} - 1)  \,,&& 
\label{V3brequalDelta3low1}      \\
t^{b/2}(t^{2k} + 1)(t^{k} + 1) \Delta_{\widehat{\gamma}}(t)  =
(t^{b}+\epsilon_b t^{k} ) (t^{2k} +t^{k} + 1 )t/G 
+ t^{b+3k} +\epsilon_b  \,,&& 
\label{V3brequalDelta3low}      \\
%%%%%%
t^{a/2}(t^{2k} + 1)(t^{k} + 1)\Delta_{\widehat{\beta}}(t) =
(  t^{a+k} + \epsilon_a )  (t^{2k} +t^{k} + 1 )t/G
+  t^{a} +\epsilon_a t^{3k}\,.&&
\label{V3brequalDelta3}
\end{eqnarray}

When $\mu(\widehat \beta)=3$ we must have $(t-1)^2 | \Delta_{\widehat{\beta}}$, p.~6 \cite{105},
contradicting (\ref{V3brequalDelta3}).

Substitute $k = 1$ in (\ref{V3brequalDelta3low}, \ref{V3brequalDelta3}) to see
$\Delta_{\widehat{\gamma}}(t) = \Delta_{T_{b + 1}}(t)= \Delta_{T_{a - 1}}(t)=
\Delta_{\widehat{\beta}}(t) $.

For $k \geq 2$, we cannot have $a=0,\,1$ as the choice $t=e^{\imath \pi/2k}$ in (\ref{V3brequalDelta3})
leads to a contradiction.
Assume now that $a$ is even, so $\widehat \beta$ and $\widehat \gamma$ are knots.
Eq.~\ref{V3brequalDelta3} at $t=e^{\imath \pi/k}$ and $t=e^{\imath \pi/2k}$ show 
that we must have $k=2$ and $a=2(4p+1)$.
%, for $p \geq 0$.
For $p=0$ we have $\Delta_{\widehat{\beta}}(t)=1$; use (\ref{V3brequalDelta1}) to see $\Delta_{\widehat{\gamma}}(t)=1$.
For $p \geq 1$, set $k=2$ and $b=8p-2$ in (\ref{V3brequalDelta3low1}) to derive
the following. 
Eq.~\ref{HJones3brDeltafor3brknotlow} is equivalent to (\ref{HJones3brDeltafor3brknotlowequiv}).

\begin{eqnarray}
t^{b/2}(t^{8} -1) \Delta_{\widehat{\gamma}}
%&=& [t^4(t^3-t^2+1)-t^3 + t^{2} -1 ](t^{8p-4} + 1) 
%+  (t^{2} -1)t^{8p-4} ( t^{8}-1)  \\
%
= (t^4-1)(t^3-t^2+1)(t^{8p-4} + 1) 
+  (t^{2} -1)t^{8p-4} ( t^{8}-1)\,.&&
\label{HJones3brDeltafor3brknotlowequiv}
\end{eqnarray}

Similarly, for odd $a$ we have $a=k(4p+3)$, for $p \geq 0$ and $k$ odd.
Set $b=4kp+k$ in (\ref{V3brequalDelta3low1}) to derive
the following. 
Eq.~\ref{HJones3brDeltafor3br2linklow} is equivalent to (\ref{HJones3brDeltafor3br2linklowequiv}).
\begin{eqnarray}
t^{b/2}(t^{4k} - 1) \Delta_{\widehat{\gamma}}(t)
&=& (t^{4kp}- 1 ) \{(t^{k} - 1) + t^{k}(t^2-t)(t^{3k}- 1 )/(t^3-1)   \} \nonumber \\
&&+ (t^{4k} - 1)t^{4kp} (t^{k} - 1)  \,.   
\label{HJones3brDeltafor3br2linklowequiv}
\end{eqnarray}
\end{proof}

\begin{prop} 
If $\beta, \gamma \in B_3$ have $V_{\widehat \beta} = V_{\widehat \gamma}$, we have 
$P_{\widehat \beta} = P_{\widehat \gamma}$\,.
\label{VequivP3br}
\end{prop}

\begin{proof}
As in Lemma~\ref{HJones3br}, write $a=w(\beta)$ and $b=w(\gamma)$ with $a=b+2k$.
By Prop.~\ref{Writhe3brPlusAny} we may assume $a > b$. 
We may also assume $a \geq 0$ (indeed $a > 0$), since
$\Delta_{\overline{L}}(t) =\Delta_{L}(1/t) $ and $V_{\overline{L}}(t) =V_{L}(1/t) $.
Lemma~\ref{HJones3br} and Prop.~\ref{HConway3br} show the result holds when $k=1$ 
or $a=k=2$.
It suffices to show that  
Eqs.~\ref{HJones3brDeltafor3brknotlow}, 
\ref{HJones3brDeltafor3br2linklow} are not satisfied by three-braid links,
with $b \geq 3$ by Lemma~\ref{HJones3br}.

%Using the expression from Prop.~\ref{Murasugi3brpartitioning} 
%\ref{Murasugi3brOmega6}, 
Eqs.~\ref{HJones3brDeltafor3brknotlow} and \ref{HJones3brDeltafor3br2linklow} show $\deg t^{b/2} \Delta_{\widehat \gamma} = b$,
so (\ref{MaxDegForDeltaInOmega0to5}) implies $\gamma \in \Omega_6^*$, say $\gamma = \alpha_2^{3d} \eta$
(see Prop.~\ref{Murasugi3brpartitioning}). 
When $d=0$, all coefficients of $\Delta_{\gamma}$ are non-zero and alternate in sign, by 
Prop.~4.2, p.~13 \cite{70}. 
%This is also
%true for $\gamma=\alpha_2^{\pm 3} \sigma_2^{-e_1} \sigma_2^{E_1}$
%when $\min(e_1,E_1)=1$, by Eq.~\ref{DeltaAnyPowerCenter}. 
Eqs.~\ref{HJones3brDeltafor3brknotlow} and \ref{HJones3brDeltafor3br2linklow} do not have these properties.

Use (\ref{DeltaAnyPowerCenter}, \ref{Deltatorusstdform}) to see
$t^{b/2} G \Delta_{\widehat {\gamma} } =
t\{   t^{1-e +6d} +\epsilon_{e+1} t^{3d-e} A_{2+e}       - \epsilon_{e} \}$ 
for $\gamma=\alpha_2^{3d} \sigma_2^{-e} \sigma_1^{1}$.
As $b=1-e +6d > 0$ and $\deg t^{b/2} \Delta_{\widehat \gamma} = b$,
the highest degree term must be $t^{1+3d-e} A_{2+e}$, so $3d= b$ and $e=3d+1$.
This gives us $t^{b/2} G \Delta_{\widehat {\gamma} } =    t^{3d+1} +\epsilon_{d} A_{3d+3} + \epsilon_{d} t$.
However, (\ref{HJones3brDeltafor3brknotlow}, \ref{HJones3brDeltafor3br2linklow}) shows that the lowest degree terms of 
$t^{b/2} G \Delta_{\widehat {\gamma} } $ are $1+t+0t^2$ and $-1-t-t^2$, respectively.
As $t^{3d+1} +\epsilon_{d} A_{3d+3} + \epsilon_{d} t$ does not have either of these patterns, 
$\gamma=\alpha_2^{3d} \sigma_2^{-e} \sigma_1^{1}$ is not allowed.
Similarly, $t^{b/2} G \Delta_{\widehat {\gamma} } =  t^{3d} A_{E+2}   +  t^{E+6d}   - \epsilon_{E}t$
for $\gamma=\alpha_2^{3d} \sigma_2^{-1} \sigma_1^{E}$.
As $b=E-1 +6d > 0$, the first term, $t^{3d} A_{E+2}$, must be of highest degree. This implies $d=0$ that is disallowed.
%we need $b=3d+E+1$ to achieve $\deg t^{b/2} \Delta_{\widehat \gamma} = b$, but
%these two values for $b$ require $3d=2$, an impossibility. 
For $\gamma = \alpha_2^{3d} \eta$
in Prop.~\ref{Murasugi3brpartitioning} \ref{Murasugi3brOmega6}, we now know $r>2$; or
$r=1$ and $e_1,\,E_1 \geq 2$. In either case we have $\sum_{k=1}^{r} E_k \geq 2$ and  $\sum_{k=1}^{r} e_k  \geq 2$.

By (\ref{Alexander3brApprox}) we have
$\deg t^{b/2} \Delta_{\widehat \eta} = (b+|\eta|-2)/2$.
Write $ t^{b/2} \Delta_{\widehat \gamma}$ as $t^{b/2} \Delta_{\widehat \eta} + F$, with
$F= t (t^{3d}-1) (t^{3d+w(\eta)}-\epsilon_{w(\eta)}) /G$.
%$F=t^{b/2} t (t^{3d}-1) (t^{3d+w(\eta)}-\epsilon_{w(\eta)}) /(Gt^{3d} t^{w(\eta)/2})$.
When $d < 0$, we have $3d+w(\eta)>b>0$, so
$\deg  F = -1 + 3d+w(\eta) $, and this exceeds $b$. 
This implies $\deg  F = \deg t^{b/2} \Delta_{\widehat \eta}$, so
$-1 + 3d+w(\eta) = (b+|\eta|-2)/2$, and thus 
$w(\eta) = |\eta|$. This is impossible for $\eta$, an alternating braid, so $d > 0$.

Write $E= \sum_{k=1}^{r} E_k$ and $e= \sum_{k=1}^{r} e_k$, so $E,\,e \geq 2$ and $d>0$. 
When $3d+w(\eta) \geq 0$, we see $\deg  F$ is below $b$, so 
$t^{b/2} \Delta_{\widehat \eta}$ contributes the highest degree term to $t^{b/2} \Delta_{\widehat \gamma}$
and $b=(b+|\eta|-2)/2$, i.e. $b=|\eta|-2 = E+e-2$.
Now $b=E+e-2=6d+E-e$ implies $e=3d+1$, so
$3d+w(\eta) \neq 0$ (lest $E=1$). 
With $3d+w(\eta) > 0$, we see $F$ is monic and $\deg  F = b-1$. 
Note that the coefficient of $t^b$ in both (\ref{HJones3brDeltafor3brknotlow}) and (\ref{HJones3brDeltafor3br2linklow})
is $1$ and the coefficent of $t^{b-1}$ is $0$. 
In order that $t^{b/2} \Delta_{\widehat \gamma}$ be monic, $e$ must be odd (see Eq.~\ref{Alexander3brApprox}).
Eq.~\ref{Alexander3brApprox} now also tells us that the $t^{b-1}$ coefficient for
$t^{b/2} \Delta_{\widehat \eta} + F$ is $-r \neq 0$, so 
we cannot have $3d+w(\eta) \geq 0$.

The final case to consider is when $3d+w(\eta) < 0$. We use the same notation as in the prior
paragraphs.
We now have $\deg  F = 3d-1$. 
If $\deg  F > b$, we must have $\deg  F = \deg t^{b/2} \Delta_{\widehat \eta}$, i.e.
$3d-1=(b+|\eta|-2)/2 $, so $E=0$, that is disallowed.
If $\deg  F < b$, we must have $\deg t^{b/2} \Delta_{\widehat \eta}=b$, so
$(b+|\eta|-2)/2 = b$, i.e. $b=|\eta|-2 = E+e-2$.
As $3d-1 < b= 6d+E-e$, it follows that $e< 3d+E+1$.
The assumption $3d+w(\eta) < 0$ means that $3d+E < e$, a contradiction, so 
we must have $\deg  F = b$.
If $w(\eta)$ is even, the coefficient of the highest power term of $F$ is $-1$.
By Eq.~\ref{Alexander3brApprox}, the coefficient of $t^b$ in $t^{b/2} \Delta_{\widehat \eta}$ 
is $0, \pm 1$, so the lead coefficient of $ t^{b/2} \Delta_{\widehat \gamma}$ is one of $0,-1,-2$.
This conflicts with (\ref{HJones3brDeltafor3brknotlow}) and (\ref{HJones3brDeltafor3br2linklow}) that
have a coefficient of $1$ for $t^b$.
Thus $w(\eta)$ is odd, so $F$ is monic and we must have 
$\deg t^{b/2} \Delta_{\widehat \eta} = (b+|\eta|-2)/2 < b = 3d-1$.
As $b+|\eta|-2 < 6d-2$, we have $6d+ w(\eta) +|\eta| <6d$, an impossibility.
\end{proof}

%%%%%%%%%%%%%%%%%%%%%%%%%%%%%%%%%%%%%  PROOFS  %%%%%%%%%%%%%%%%%%%%%%%%%%%%%%%%%%%%%%
%%%%%%%%%%%%%%%%%%%%%%%%%%%%%%%%%%%%%  PROOFS  %%%%%%%%%%%%%%%%%%%%%%%%%%%%%%%%%%%%%%
\appendix
\section{Expression for Laurent MFW polynomials}
This section contains the proof of Thm.~\ref{Hcoeffrel} whose inclusion in Section~\ref{SectionMain}
would provide no important insight and would 
interrupt the flow.

%%%%%%%%%%%%%%%%%%%%%%%%%%  Expression for Laurent coefficient polynomials %%%%%%%%%%%%%%%%%%%%%%%%%%%%%%%%%%

\noindent
\textbf{Theorem~\ref{Hcoeffrel}.}
%\begin{thm}
For an arbitrary braid, $\beta$, of $n \geq 1$ strands, the \mbox{Homflypt} polynomial for 
$\widehat{\beta}$ is given by (\ref{HLaurentform}) with:
\begin{enumerate}%[(i)]
\item $h_0 = z^{-2}\{\,C_{w+4-n} - \nabla_{\widehat{\beta}}\} -  q_{0} $\,,
with $q_{0} = \sum_{j=3}^{n-1}  \,z^{-2}(C_{2j-3}-1) \,h_j$\,,
%\item $p_0 = z^{n-3}\{\,C_{w+4-n} - \nabla_{\widehat{\beta}}\} - q_{0} $\,,
%with $q_{0} = \sum_{j=3}^{n-1}  \,z^{-2}(C_{2j-3}-1) \,p_j$\,,
%
\item $h_{1} = \nabla_{\widehat{\beta}} -h_0 - h_2 - q_{1} $\,,
with $q_{1} = \sum_{j=3}^{n-1} \,h_j$\,,
%\item $p_{1} = z^{n-1} \,\nabla_{\widehat{\beta}} -p_0 - p_2 - q_{1} $\,,
%with $q_{1} = \sum_{j=3}^{n-1} \,p_j$\,,
%
\item $h_2 = z^{-2}\{\,C_{w+2-n} - \nabla_{\widehat{\beta}}\}  - q_{2}  $\,,
with $q_{2} = \sum_{j=3}^{n-1}  \,\,z^{-2}(C_{2j-1}-1) \,h_j$\,.
%\item $p_2 = z^{n-3}\{\,C_{w+2-n} - \nabla_{\widehat{\beta}}(z)\}  - q_{2}  $\,,
%with $q_{2} = \sum_{j=3}^{n-1}  \,\,z^{-2}(C_{2j-1}-1) \,p_j$\,.
%
\end{enumerate}
%\end{thm}
\label{ProofHcoeffrel}

\begin{proof}
The proof uses similar relations as in Remark 2, p.~410, \cite{69} that is the proof for the
Corollary there ($\sqrt{t}$ is replaced by $s$ below). 
A slight modification is needed as the skein relation defining the 
\mbox{Homflypt} polynomial in \cite{69} is $a^{-1} P_{L_+}(a,z) -a P_{L_-}(a,z) +z P_{L_0}(a,z) = 0$,
and this differs from (\ref{Hskeinrel}).

Let us start with (\ref{HLaurentform}) and choose $z=s-s^{-1}$.
When $v=1$, we obtain the Alexander polynomial, $\Delta_{\widehat{\beta}}(s^2)$.
Appeal now to Thm.~8.4.1 parts 5 and 6 (p.~108, \cite{54}), and the fact that
$w-n \equiv \mu(\widehat{\beta}) \mbox{ mod } 2$.
When $v=s$, we obtain $\epsilon_{w-n+1}$. 
Similarly, when $v=s^{-1}$, we obtain $1$.
For these values, with $n \geq 3$, Eq.~\ref{HLaurentform} may be represented in matix form as:
\begin{eqnarray*}
\left ( 
\begin{array} {cc} \Delta_{\widehat{\beta}}(s^2) \\ \epsilon_{w-n+1} \\ 1
\end{array}
\right ) &=&
\left [ 
\begin{array} {cccc} 1          & 1          & \cdots & 1  \\
                     s^{w-n+1}  & s^{w-n+3}  & \cdots & s^{w+n-1} \\
                     s^{-w+n-1} & s^{-w+n-3} & \cdots & s^{-w-n+1}
\end{array}
\right ]
\left ( 
\begin{array} {cc} h_0(z) \\ \cdots \\  h_{n-1}(z)
\end{array}
\right ) .
\end{eqnarray*}

If we designate by $A$ the $3 \times 3$ matrix consisting of the first three columns,
and by $B$ the matrix consisting of the remaining $n-3$ columns we have 
\begin{eqnarray*}
\left ( 
\begin{array} {cc} \Delta_{\widehat{\beta}}(s^2) \\ \epsilon_{w-n+1} \\ 1
\end{array}
\right ) &=& A 
\left ( 
\begin{array} {cc} h_0(s-s^{-1}) \\ h_1(s-s^{-1})\\  h_{2}(s-s^{-1})
\end{array}
\right ) + B
\left ( 
\begin{array} {cc} h_3(s-s^{-1}) \\ \cdots \\  h_{n-1}(s-s^{-1})
\end{array}
\right ) .
\end{eqnarray*}

We may rewrite this as 
\begin{eqnarray*}
\left ( 
\begin{array} {cc} h_0(s-s^{-1}) \\ h_1(s-s^{-1}) \\ h_2(s-s^{-1})
\end{array}
\right ) &=&
A^{-1}
\left ( 
\begin{array} {cc} \Delta_{\widehat{\beta}}(s^2) \\ \epsilon_{w-n+1} \\ 1
\end{array}
\right ) - 
A^{-1} B
\left ( 
\begin{array} {cc} h_3(s-s^{-1}) \\ \cdots \\  h_{n-1}(s-s^{-1})
\end{array}
\right ) .
\end{eqnarray*}

It is readily seen that with $d(s)= (s^4-1)(s^2-1)$ we have
\begin{eqnarray*}
d(s) A^{-1} &=& 
\left [ 
\begin{array} {ccc}  -s^2(1+s^2)    & s^{-w+n-1}         & s^{w-n+7}         \\
                     (1+s^2)(1+s^4) & -s^{-w+n-1}(1+s^2) & -s^{w-n+5}(1+s^2) \\
                     -s^2(1+s^2)    & s^{-w+n+1}         &s^{w-n+5} 
\end{array}
\right ].
\end{eqnarray*}

This implies
\begin{eqnarray*}
d(s) h_0(s-s^{-1}) &=& -s^2(1+s^2) \Delta_{\widehat{\beta}}(s^2) + s^{-w+n-1} \epsilon_{w-n+1} + s^{w-n+7}  \\
&&- \sum_{j=3}^{n-1} (s^{2j}-s^4-s^2+s^{6-2j})h_j\,,   \\
d(s) h_1(s-s^{-1}) &=& (1+s^2)(1+s^4) \Delta_{\widehat{\beta}}(s^2)  -s^{-w+n-1}(1+s^2) \epsilon_{w-n+1} \\  
&& -s^{w-n+5}(1+s^2) + \sum_{j=3}^{n-1} (1+s^2)(s^{2j}-s^4-1+s^{4-2j})h_j\,,  \\
d(s) h_2(s-s^{-1}) &=& -s^2(1+s^2) \Delta_{\widehat{\beta}}(s^2) + s^{-w+n+1}  \epsilon_{w-n+1} + s^{w-n+5} \\
&&- \sum_{j=3}^{n-1} (s^{2+2j}-s^4-s^2+s^{4-2j})h_j\,.
\end{eqnarray*}

This is exactly the result implied by Thm.~\ref{Hcoeffrel}. 
As $z=s - s^{-1}$ is invertible for $|s| >1$ in the complex plane, $P_{\widehat{\beta}}(v,z)$ has the form claimed.
\end{proof}

\section*{Acknowledgements}
The author would like to thank prior referees for additional citations and possible techniques to 
simplify the proof of Theorem~\ref{Hcoeffrel}, one of which is incorporated above.  

%%%target range end to move body

%% The Appendices part is started with the command \appendix;
%% appendix sections are then done as normal sections
%% \appendix

%% \section{}
%% \label{}

%% References
%%
%% Following citation commands can be used in the body text:
%% Usage of \cite is as follows:
%%   \cite{key}          ==>>  [#]
%%   \cite[chap. 2]{key} ==>>  [#, chap. 2]
%%   \citet{key}         ==>>  Author [#]

%% References with bibTeX database:

\bibliographystyle{model1-num-names}
% 7/25/2011 suppress use of bib database below
%\bibliography{<your-bib-database>}

\begin{thebibliography}{00}

%% \bibitem must have the following form:
%%   \bibitem{key}...
%%

% \bibitem{}

\bibitem{5} J.S.~Birman, On the Jones Polynomial of closed 3-braids,
{\it Invent.  Math.} {\bf 81(2)}, (1985), 287--294.

\bibitem{10} J.S.~Birman and T.E.~Brendle, Braids: A  
Survey, 
{\it arXiv:math.GT/0409205v2} (2005).

%\bibitem{21} J.C.~Cha and C.~Livingston, KnotInfo: Table of Knot Invariants,
%{\it http://www.indiana.edu/~knotinfo}, 2010

\bibitem{30} J.~Franks and R.F.~Williams, Braids and the 
Jones polynomial, 
{\it Trans. Amer. Math. Soc.}, {\bf 303} (1987), 97--108.


\bibitem{41} V.F.R.~Jones, Hecke Algebra Representations of 
Braid Groups and Link Polynomials
{\it Ann. of Math.}, {\bf 126} (1987), 335--388.

\bibitem{51} T.~Kanenobu, An Evaluation of the Coefficient Polynomials
of the HOMFLY Polynomial of a Link,
in {\it Knots in Hellas '98, Proc. of the Int. Conf. on
Knot Theory}, Eds. C. McA. Gordon et al, World Scientific 2000,
pp.~131-137.

\bibitem{53} A.~Kawauchi, On Coefficient Polynomials of the Skein
Polynomial of an Oriented Link,
{\it Kobe J. Math.}, {\bf 11}
(1994), 49--68.

\bibitem{54} A.~Kawauchi, {\it A Survey of Knot Theory} 
(Birkhauser, Boston, 1996).

\bibitem{60} H.R.~Morton, Seifert circles and knot polynomials,
{\it Math. Proc. Cambridge Philos. Soc.}, {\bf 99}
(1986), 107--109.

\bibitem{69} H.~Murakami, A formula for the two-variable link polynomial. 
{\it Topology}, {\bf 26(1)} (1987), 409--412.

\bibitem{70} K.~Murasugi, On Closed 3-braids
{\it Mem. of Amer. Math. Soc.}, {\bf 151} (1974), 
American Mathematical Society, Providence, RI.

\bibitem{71} K.~Murasugi, {\it Knot Theory \& Its 
Applications} (Birkhauser, Boston, 1996).

\bibitem{105} A.~Stoimenow, Properties of Closed  
3-Braids and Other Link Braid Representations
{\it arXiv:math.GT/0606435.v2} (2007).

% 7/25/2011 enable non-bib database
\end{thebibliography}

%% Authors are advised to submit their bibtex database files. They are
%% requested to list a bibtex style file in the manuscript if they do
%% not want to use model1-num-names.bst.

%% References without bibTeX database:
% 7/25/2011 enable non-bib database

\end{document}